\documentclass[submit]{aiaa-tc}
\pdfoutput=1

 \usepackage{amsmath} 
 \usepackage{amssymb} 
 \usepackage{mathrsfs}
 \usepackage{multicol}
 \usepackage{mathtools}

 \AIAAcopyright{\AIAAcopyrightD{2017}}

\bibdata

\usepackage{amssymb}
\usepackage{graphicx}
\usepackage{mathtools}
\usepackage[makeroom]{cancel}
\usepackage{amsthm}
\usepackage{float}
\usepackage{subfig}
\graphicspath{ {./} }

\newcommand{\va}{\boldsymbol{a}}
\newcommand{\bphi}{\boldsymbol{\phi}}
\newcommand{\vx}{\boldsymbol{x}}
\newcommand{\vp}{\boldsymbol{p}}
\newcommand{\mM}{\boldsymbol{M}}

\newcommand{\vr}{\boldsymbol{r}}
\newcommand{\mC}{\boldsymbol{C}}
\newcommand{\mW}{\boldsymbol{W}}
\newcommand{\hmW}{\hat{\boldsymbol{W}}}
\newcommand{\vy}{\boldsymbol{y}}
\newcommand{\vz}{\boldsymbol{z}}
\newcommand{\mH}{\boldsymbol{H}}
\newcommand{\vv}{\boldsymbol{v}}
\newcommand{\mI}{\boldsymbol{I}}
\newcommand{\hLambda}{\hat{\Lambda}}
\newcommand{\Nboot}{N_{\text{boot}}}
\newcommand{\hvw}{\hat{\boldsymbol{w}}}
\newcommand{\bmat}[1]{\begin{bmatrix}#1\end{bmatrix}}
\newcommand{\mini}[1]{\underset{#1}{\mathrm{minimum}}}

\begin{document}

\title{Active Subspaces of Airfoil Shape Parameterizations}


 \author{
  Zachary J. Grey%
    \footnote{Applied Mathematics and Statistics Department, Chauvenet Hall 1015 14th St.}
    \footnote{Rolls-Royce Corporation, Engineering Sciences -- Robust Design Technical Lead}
  \ and Paul G. Constantine\thanksibid{1}\\
  {\normalsize\itshape
   Colorado School of Mines, Golden, Colorado, 80401, USA}\\
 }
\maketitle

\begin{abstract}
Design and optimization benefit from understanding the dependence of a quantity of interest (e.g., a design objective or constraint function) on the design variables. A low-dimensional \emph{active subspace}, when present, identifies important directions in the space of design variables; perturbing a design along the active subspace associated with a particular quantity of interest changes that quantity more, on average, than perturbing the design orthogonally to the active subspace. This low-dimensional structure provides insights that characterize the dependence of quantities of interest on design variables. Airfoil design in a transonic flow field with a parameterized geometry is a popular test problem for design methodologies. We examine two particular airfoil shape parameterizations, PARSEC and CST, and study the active subspaces present in two common design quantities of interest, transonic lift and drag coefficients, under each shape parameterization. We mathematically relate the two parameterizations with a common polynomial series. The active subspaces enable low-dimensional approximations of lift and drag that relate to physical airfoil properties. In particular, we obtain and interpret a two-dimensional approximation of both transonic lift and drag, and we show how these approximation inform a multi-objective design problem.
\end{abstract}

\section*{Nomenclature}
\noindent\begin{tabular}{@{}lcl@{}}
$\mathcal{A}$ &=& airfoil subset\\
$\partial \mathcal{A}$ &=& airfoil surface/boundary\\
$\boldsymbol{x}$  &=& parameter vector \\
$\boldsymbol{y}$  &=& active variable vector \\
$\boldsymbol{z}$  &=& inactive variable vector\\
$\mW_1$  &=& orthonormal basis of the active subspace \\
$\hmW_1$  &=& approximation of $\mW_1$ \\
$\hat{\boldsymbol{w}}_i$  &=& $i^{th}$ column of $\hmW_1$  \\
$f$ &=& quantity of interest \\
$s$ &=& shape function \\
$c_f$ &=& CST class functions \\
$s_f$ &=& CST shape functions \\
$\mM$ &=& PARSEC constraint matrix\\
$\boldsymbol{p}$ &=& PARSEC constraint vector \\
$\boldsymbol{\phi}$ &=& basis function vector \\
$s_U$ &=& upper shape function \\
$s_L$ &=& lower shape function \\
$\ell$   &=& scaled length coordinate \\
$t$   &=& intermediate length scale \\
$\boldsymbol{a}, \boldsymbol{b}$   &=& shape coefficient vectors \\
\end{tabular} \\

\section{Introduction}

The approximation of transonic flow solutions modeled by Euler's equations has remained a relevant test problem in aerospace engineering---specifically, the design and evaluation of airfoils in a two-dimensional inviscid compressible flow field. The numerical representation of airfoils in transonic flow regimes involves parametrized geometries that inform optimization procedures, as described by Vanderplaats\cite{Vanderplaats2007}. Several applications---e.g., Garzon\cite{Garzon2003}, Dow and Wang\cite{Dow2014}, and Lukaczyk et al.\cite{Lukaczyk2014}---approximate flow field solutions with parametric representations of airfoil shapes to produce estimates of quantities of interest (QoI). The performance of a design is measured by studying changes in QoI (e.g., dimensionless force quantities like lift and drag or measures of efficiency) as shape parameters vary. The QoI are computed from sufficiently accurate flow field approximations. The flow field approximations depend on the parameterized airfoil, which defines the boundaries in the governing equations. Samareh\cite{Samareh2001} surveys parameterization techniques for high fidelity shape optimization. H\'ajek\cite{Hajek} also offers a concise taxonomy of parameterizations and optimization routines.

We seek to further understand the QoI dependence on airfoil shape parameters through transformations and approximations that inform algorithms, insights, and pedigree---and ultimately improve airfoil designs. Toward this end, we search for \emph{active subspaces} in the map from airfoil shape parameters to design quantities of interest, lift and drag, in a two-dimensional transonic airfoil model. An active subspace is defined by a set of important directions in the space of input parameters; perturbing parameters along these important directions (when present) changes the outputs more, on average, than perturbing the parameters along unimportant directions (i.e., directions orthogonal to the important directions). Such structure, when present, can be exploited to reduce the input space dimension and subsequently reduce the cost of algorithmic approaches to the design problem. Our work builds on recent work exploiting active subspaces in shape optimization\cite{Lukaczyk2014}, vehicle design\cite{othmer2016}, jet nozzle design\cite{SEQUOIA2017}, and dimensional analysis\cite{delRosario2017}, as well as a gradient-based dimension reduction approach from Berguin and Mavris\cite{Berguin2015}. 

The present work extends previous work as follows. We study two common shape parameterizations---PARSEC\cite{Sobieczky1999.} and CST\cite{Kulfan2008}---and we show how design-enabling insights derived from active subspaces are consistent across the two parameterizations. Moreover, we show how the active subspace reduces the dimension and enables simple visualization tools that lead to deep insights into the lift/drag Pareto front in a multi-objective optimization for the airfoil design problem.

\section{Airfoil Design}
\subsection{An interpretation of airfoils} \label{sec:airfoil}
An airfoil is a subset of a two-dimensional Euclidean space,
$$
\mathcal{A} \subset \mathbb{R}^2,
$$
with a boundary, $\partial \mathcal{A}$, representing the surface of the airfoil. The airfoil surface can be partitioned into two shape functions, $s_U$ and $s_L$, representing the upper and lower surfaces, respectively. Each shape function maps a normalized spatial coordinate $\ell\in[0,1]$, running from leading to trailing edge, to its respective surface. The following are necessary conditions on $s_U(\ell)$ and $s_L(\ell)$ such that they represent an \emph{airfoil} for the present study:
\begin{enumerate} \label{props}
\item \emph{Bounded}: $0 \leq s_U(\ell)\leq U<\infty$ for an upper bound $U$ and $-\infty < L \leq s_L(\ell) \leq 0$ for a lower bound $L$.
\item \emph{Smooth}: continuous and twice differentiable for $\ell\in(0,1)$.
\item \emph{Fixed endpoints}: $s_U(0) = s_L(0) = 0$ and, optionally, $s_U(1) = s_L(1)= 0$.
\item \emph{Feasible}: $s_U(\ell) > s_L(\ell)$ for $\ell \in (0,1)$.
\end{enumerate}
These conditions can be interpreted physically. Boundedness ensures that the resulting airfoil has a minimum and a maximum, i.e., the shape must remain finite to be manufacturable. Smoothness ensures that the interior of the shape is manufacturable to manufacturing precision. Traditionally, airfoils have a round leading edge, which produces an unbounded first derivative in $s_U$ and $s_L$ at $\ell=0$. The fixed endpoints ensure the two surfaces meet at $\ell=0$ and $\ell=1$. A sharp trailing edge is not manufacturable to finite precision, though it is used in most computational models. Feasibility ensures surfaces do not intersect. We ensure property four is satisfied by assuming positive shapes, $s_U(\ell)>0$ and $-s_L(\ell)>0$ for $\ell\in(0,1)$. 

\subsection{Airfoils for design}
\begin{figure}[H] 
\centering
\subfloat{
\includegraphics[width=0.55\textwidth]{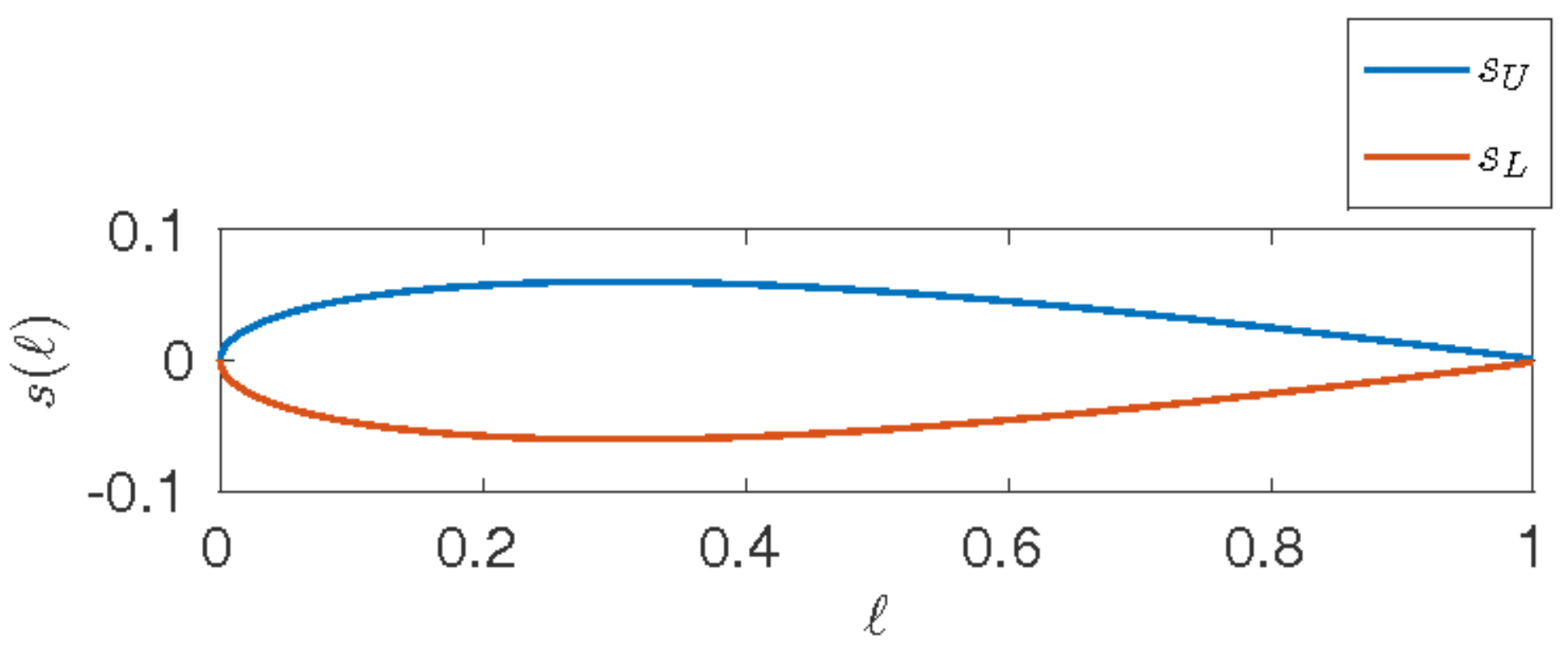}
}
\caption{NACA 0012 airfoil represented as two shapes: the upper surface $s_U(\ell)$ and the lower surface $s_L(\ell)$. }
\label{fig:NACA0012}
\end{figure}

Airfoil shapes can be parameterized by (i) deformation of a nominal shape\cite{Wu2003}, (ii) smooth perturbations to a nominal shape\cite{Hicks1978}, (iii) Karhunen-Lo\`eve expansions\cite{Garzon2003}, (iv) Joukowski\cite{Hajek} or Theodorsen-Garrick conformal mappings\cite{Iollo1999}, and (v) piece-wise spline interpolation with relevant bounding\cite{Song2011}. Alternatively, engineering intuition has motivated a variety of parameterized shapes, specifically those from the National Advisory Committee for Aeronautics (NACA)\cite{Jacobs1933}, Sobieczky\cite{Sobieczky1999.}, and Kulfan\cite{Kulfan2008}. Considering the engineering insights, we study partitions of $\partial \mathcal{A}$ into shape functions that can be written as a linear combination of basis functions, $\phi_j:\mathbb{R} \rightarrow \mathbb{R}$, such that, for $s$ as one of $s_U$ or $s_L$,
\begin{equation} \label{genshape}
s(\ell; \{a_j\}) = \sum_{j=1}^{k} a_j\phi_j(\ell) \quad \text{where} \quad a_j \in \mathbb{R}, \quad \ell \in [0,1].
\end{equation}
Varying the coefficients $a_j$ modifies the shape functions. Hence, the coefficients $\va = [a_1,\dots,a_k]^T$ define the shapes with basis $\bphi(\ell)=[\phi_1(\ell),\dots,\phi_k(\ell)]^T$.
Vector notation enables a concise expression for separate expansions in \eqref{genshape}. For example, the top surface of a NACA airfoil can be expressed as\cite{Jacobs1933}
\begin{equation} \label{NACA}
s_U(\ell) = \tau \,\va^T\bphi(\ell), \quad \phi_j(\ell) = \begin{cases} \sqrt{\ell}, & j=1 \\ \ell^{j-1}, & j=2,\dots,5\end{cases}, \quad \tau > 0.
\end{equation}
For NACA airfoils, the coefficients $\va$ sum to zero, which implies the fixed endpoint at the trailing edge. And $\va$ are constrained to ensure feasibility ($s_U>s_L$). The basis expansion is smooth on the open interval, bounded on the closed interval, and admits a leading edge ($\ell=0$) derivative singularity to model the roundness. Figure \ref{fig:NACA0012} shows the well-known NACA 0012 as two shape functions. 

\begin{figure}[H]
\centering
\subfloat[Round nose graph in $\ell$]{
\includegraphics[width=0.45\textwidth]{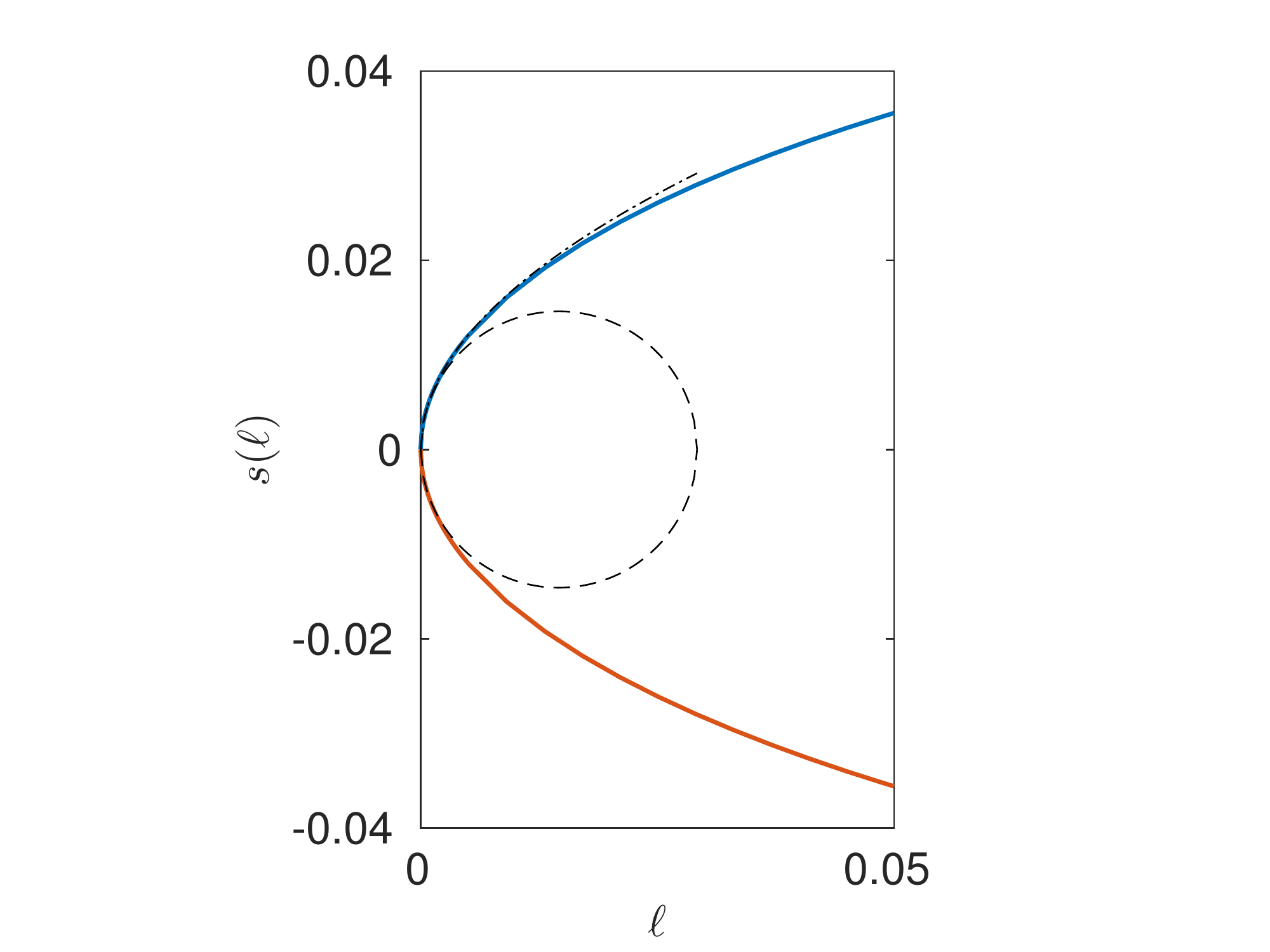}
}
\hfil
\subfloat[Round nose graph in $t=\sqrt{\ell}$]{
\includegraphics[width=0.45\textwidth]{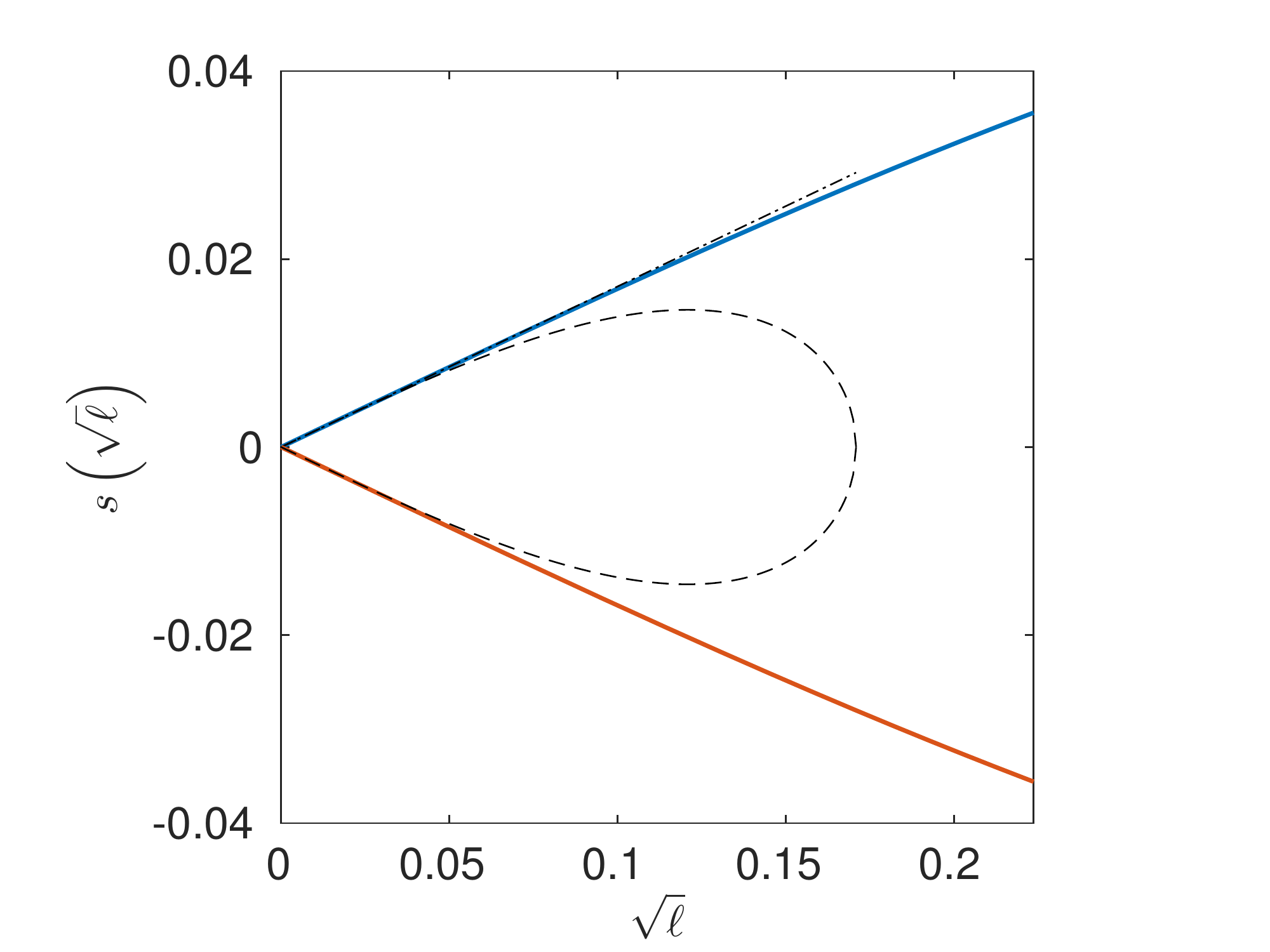}
}
\caption{Approximating the round nose of the airfoil using the original spatial coordinate $\ell$ and the change-of-variables $t=\sqrt{\ell}$.}
\label{fig:round}
\end{figure}

Consider $0<\epsilon \ll 1$ and $\ell \in [0,\epsilon)$ such that powers of $\ell$ are negligibly small. In this neighborhood of $\ell=0$, the upper surface can be approximated by
\begin{equation}
s_U(\ell) \approx \tau a_1\sqrt{\ell}  \quad \ell\in [0,\epsilon).
\end{equation}
That is, the symmetric NACA airfoil behaves like a multiple of the square root and admits the round-nose leading edge (consequently, a singularity at $\ell=0$). Figure \ref{fig:round} shows this behavior along with a $t=\sqrt{\ell}$ change-of-variables for perspective. Writing the NACA 0012 basis functions $\phi_j(\ell)$ from \eqref{NACA} with the change-of-variables, 
\begin{equation}\label{NACAt}
\phi_j(t) = \begin{cases} t, & j=1 \\ t^{2(j-1)}, & j=2,\dots,5\end{cases}.
\end{equation}
The transformation from $t$ to $\ell$ helps compare shape representations in the following sections.

Figure \ref{fig:round} demonstrates the utility of approximating shapes with the $\sqrt{\ell}$ change of variables. This transformation also appears in other shape function representations that have circular leading edge (i.e., round-nose) airfoils: Parametric Sections (PARSEC)\cite{Sobieczky1999.} and Class-Shape Transformations (CST)\cite{Kulfan2008}, which we discuss in the following sections. These parametric representations admit polynomial series in $\sqrt{\ell}$ that achieve the round-nose shapes. 
Under the change of variables $t=\sqrt{\ell}$, the derivative singularities at $\ell=0$ for such series is verified by the chain rule,
\begin{equation}
\frac{d}{d\ell}s(\ell) = \frac{dt}{d\ell}\frac{d}{dt}s(t) \implies \frac{d}{d\ell}s(\ell) = \frac{1}{2\sqrt{\ell}}\frac{d}{dt}s(t).
\end{equation}

It is possible to eliminate the singularity in this shape partition through conformal mapping or additional partitioning of $\partial \mathcal{A}$. A classical representation involves the Joukowski conformal mapping of a circle but the parameterization admits a limited space of shapes based on two parameter values\cite{Abbott1959}. Designing inverse conformal mapping parameterizations can be difficult and existing conformal mappings are often associated with shapes admitting undesirable aerodynamic features\cite{Melin2011}. However, the Theodorsen-Garrick conformal mapping offers a more generic shape representation from the inverse of near-circles\cite{schinzinger2003,Abbott1959}. 

\subsection{Parametric shape functions}
Consider the following parameterization of the NACA airfoil, where two parameters, $x_1>0$ and $x_2>0$, independently scale the coefficients of the series \eqref{NACA},
\begin{equation} \label{NACAparm}
s_U(\ell) = (x_1\,\va)^T\bphi(\ell), \quad
s_L(\ell) = -(x_2\,\va)^T\bphi(\ell).
\end{equation}
We think of \eqref{NACAparm} as a map from a two-dimensional parameter domain to airfoils, where the parameters $\vx=[x_1,x_2]^T$ control the thickness. Using this parameterization, we can model asymmetric shapes 
with different values for $\vx$ (that is, $x_1\not=x_2$). Given a shape, we can compute a designer's quantity of interest from a flow field approximated on the domain containing the airfoil; we denote the quantity of interest $f$. Consider the composition map,
\begin{equation}
\label{eq:composition}
\vx \longmapsto s_U\, , \, s_L \longmapsto f,
\end{equation}
which suggests the following numerical procedure:
\begin{enumerate}
\item Choose parameters $\vx$.
\item Compute parametric shape functions, $s_U(\ell)$ and $s_L(\ell)$.
\item Create a mesh on the set $\mathcal{A}^c$, and compute the CFD flow field.
\item Approximate the quantity of interest, $f$, from the CFD flow field solution.
\end{enumerate}
Typically, asymmetric NACA airfoil parameterizations use an integer encoding for thickness and camber functions. The thickness and camber are combined using a polar coordinate transformation to add thickness orthogonal to a camber curve. The simpler NACA parameterization in \eqref{NACAparm} is meant to introduce the function, $f:\mathbb{R}^2 \rightarrow \mathbb{R}$, that maps shape parameters $\vx=[x_1,x_2]$ to the quantity of interest.

Given (i) a shape parameterization with $m$ parameters $\vx\in\mathbb{R}^m$ and (ii) a composition \eqref{eq:composition}, a design optimization seeks $\vx^\ast\in\mathbb{R}^m$ such that $f(\vx^\ast) \leq f(\vx)$ subject to constraints\cite{Vanderplaats2007}. However, if the number $m$ of parameters is too large, then optimization algorithms may be computationally prohibitive\cite{Lukaczyk2014}. If accurate approximations for all quantities of interest (i.e., objectives and constraints) can be obtained using a subspace of sufficiently small dimension (e.g., dimension less than or equal to two) in the shape parameter space, then visualization tools can aid the search for the optimum by providing insights into the relationship between inputs and outputs. In section \ref{sec:results}, we discover two-dimensional active subspaces that enable these visualizations.

\subsection{PARSEC airfoils}
The PARSEC parameterization defines a set of shape functions; this parameterization was popularized by Sobieczky\cite{Sobieczky1998.} for applications in transonic flow regimes. There have been several adaptations to PARSEC---e.g., Derksen and Rogalsky\cite{DerksenR.} adapt PARSEC for an optimization routine, while Wu et al.\cite{Wu2003} focus on comparisons for a specific optimization problem. We focus on the original PARSEC parameterization for simplicity. The PARSEC expansion for a generic shape $s(\ell)$ is
\begin{equation} \label{PARS}
s(\ell) \;=\; s(\ell; \va) \;=\; 
\sum_{j=1}^{6} a_j\ell^{j-1/2} \;=\;
\va^T\bphi(\ell),
\qquad \ell\in[0,1],
\end{equation}
with coefficients $\va=[a_1,\dots,a_6]^T$ and basis functions $\bphi(\ell)=[\ell^{1/2},\dots,\ell^{11/2}]^T$. Each surface $s_U$ and $s_L$ has its own independent coefficient vector, $\va_U$ and $\va_L$.

If we substitute the change-of-variables $t=\sqrt{\ell}$, the basis expansion becomes an odd-power polynomial in $t$,
\begin{equation} \label{PARZ}
s(t; \va) \;=\; \sum_{j=1}^{6} a_jt^{2j-1}, \qquad t \in [0,1].
\end{equation}
This change-of-variables enables precise comparison with the NACA 0012 parameterization \eqref{NACAparm} with basis \eqref{NACAt}. In particular, PARSEC can not exactly reproduce a NACA 0012 airfoil because of the difference in bases; PARSEC uses odd powers of $t$ while NACA 0012 includes even powers of $t$ with fixed, nonzero coefficients. Moreover, the change-of-variables gives insight into the ability of PARSEC to produce round-nose airfoils; see Figure \ref{fig:round}.

The coefficients of both PARSEC surfaces, upper $\va_U$ and lower $\va_L$, solve a system of 12 parameterized linear equations. Each equation represents a physically interpretable characteristic of the shape, and the parameters affect the shape in physically meaningful ways. Table \ref{tab:PARSECparms} gives a description of each of the 11 parameters, denoted $\vx=[x_1,\dots,x_{11}]^T$, that characterize the PARSEC shapes.

\begin{table}[h!]
\centering
\caption{PARSEC airfoils are characterized by eleven parameters, each affecting a physical characteristic of the airfoil.}
\begin{tabular}{l|l} 
 Parameter & Description \\ [0.5ex] 
 \hline
 $x_1$ & upper surface domain coord.  \\ 
 $x_2$ & lower surface domain coord.  \\
 $x_3$ & upper surface shape coord. \\
 $x_4$ & lower surface shape coord.  \\
 $x_5$ & trailing edge shape coord. \\ 
 $x_6$ & trailing edge half-thickness\\
 $x_7$ & wedge angle\\
 $x_8$ & wedge half-angle \\
 $x_9$ & upper surface inflection\\
 $x_{10}$ & lower surface inflection\\
 $x_{11}$ & leading edge radius\\[1ex] 
\end{tabular}
\label{tab:PARSECparms}
\end{table}

In PARSEC, bounds on the parameters $\vx$ restrict size and shape variation for design and optimization problems. Certain parameters---namely $x_5$ (trailing edge offset), $x_6$ (trailing edge half-thickness), $x_7$ (wedge angle), $x_8$ (wedge half-angle), and $x_{11}$ (leading edge radius)---affect both upper and lower surfaces. Setting $x_5=x_6=0$ produces a sharp trailing edge.The parameterization contains two nonlinear relationships: (i) trailing edge slope and (ii) leading edge radius. 

Given $\vx$, the PARSEC coefficients $\va_U=\va_U(\vx)$ and $\va_L=\va_L(\vx)$ satisfy
\begin{equation}
\label{eq:PARSECcoeff}
\begin{bmatrix} 
\mM_U(x_1) & 0 \\ 0 & \mM_L(x_2)
\end{bmatrix} 
\begin{bmatrix} 
\va_U(\vx) \\ \va_L(\vx)
\end{bmatrix} 
\;=\; 
\begin{bmatrix} 
\vp_U(\vx) \\ \vp_L(\vx) 
\end{bmatrix},
\end{equation}
where each diagonal block of the matrix is 6-by-6 and depends on either $x_1$ or $x_2$ from Table \ref{tab:PARSECparms}. All parameters from Table \ref{tab:PARSECparms}, except $x_1$ and $x_2$, affect the right hand side $[\vp_U(\vx), \vp_L(\vx)]^T$. For reasonable choices of $x_1$ and $x_2$, the matrix in \eqref{eq:PARSECcoeff} is invertible. Further details on the equality constraints are in Appendix \ref{sec:app1}. 

\subsection{Class-shape transformation airfoils} \label{sec:CST}
A more recent airfoil parameterization is the Class-Shape Transformation (CST), which was popularized by Kulfan\cite{Kulfan2008}. The shape function is decomposed into the product of a \emph{class function} $c_f$ and a smooth \emph{shape functions} $s_f$:
\begin{equation} \label{CSTCS}
s(\ell;r_1,r_2,\vx) 
= c_f(\ell; r_1,r_2)\,s_f(\ell; \vx),
\end{equation}
where $r_1$ and $r_2$ are between 0 and 1, and the vector $\vx \in \mathbb{R}^{m}$ contains shape parameters as polynomial coefficients. Each class function takes the form,
\begin{equation} \label{eq:classf}
c_f(\ell,\vr) \;=\; \ell^{r_1}\,(1-\ell)^{r_2},
\end{equation}
where $r_1$ and $r_2$ control the shape at the leading and trailing edges, respectively. For simplicity, we fix $r_{1}=1/2$ and $r_{2}=1$, which produces round nose airfoils with sharp trailing edges. Alternative class function forms have been proposed and implemented in other applications\cite{Lane2009,Powell2010,Salar2014}. The shape function is a polynomial of degree $m-1$ in $\ell$,
\begin{equation}
\label{eq:shapef}
s_f(\ell,\vx) 
\;=\; 
\sum_{j=0}^{m-1} x_j\ell^j, \quad \vx = [x_0, \dots,x_{m-1}]^T.
\end{equation}
With (i) our choices $r_1=1/2$ and $r_2=1$ and (ii) the change of variables $t=\sqrt{\ell}$, each CST surface can be written
\begin{equation}
\label{eq:cstseriest}
\begin{aligned}
s(t; \vx) &= \sum_{j=0}^{m-1} x_j\,t^{2j}\,t\,(1-t^2)\\
&= x_0 t - x_{m-1} t^{2m+1} + \sum_{j=1}^{m-1} (x_{j+1}-x_j)\,t^{2j+1}.
\end{aligned}
\end{equation}
Compare \eqref{eq:cstseriest} to the PARSEC series \eqref{PARS}; for $m=5$ in \eqref{eq:cstseriest}, both are polynomials of degree 11 in $t$ containing only odd powers of $t$. However, the parameterizations differ dramatically. The PARSEC parameters correspond to physically meaningful characteristics of the airfoil; see Appendix \ref{sec:app1}. In contrast, the parameters of CST directly impact the coefficients in the polynomial in $t$---as opposed to through a system of equations \eqref{eq:PARSECcoeff}.

The domain for the composition \eqref{eq:composition} is the space of coefficients $\vx=[\vx_U,\vx_L]^T\in\mathbb{R}^{2m}$ that parameterize the upper and lower surfaces as in \eqref{eq:shapef}. We choose $m=5$ in the numerical experiments in section \ref{sec:results} to maintain the comparison between polynomial representations of CST and PARSEC; thus, the parameter space for CST is 10-dimensional. Given a point in the 10-dimensional parameter space, the CST shape is constructed via \eqref{eq:shapef}, \eqref{eq:classf}, and \eqref{CSTCS}. Then a CFD solver estimates the flow fields and produces the quantity of interest.


\section{Active subspaces}

We closely follow the development of active subspaces as in Constantine\cite{Constantine2015}. Active subspaces are part of an emerging set of tools for dimension reduction in the input space of a function of several variables. Given (i) a probability density function $\rho:\mathbb{R}^m\rightarrow\mathbb{R}_+$ that identifies important regions of the input space and (ii) a differentiable and square-integrable function $f:\mathbb{R}^m \rightarrow \mathbb{R}$. Assume, without loss of generality, that $\rho(\vx)$ is such that
\begin{equation}
\label{eq:assump}
\int \vx\,\rho(\vx)\,d\vx = 0, \qquad
\int \vx\,\vx^T\,\rho(\vx)\,d\vx = \mI.
\end{equation}
That is, we assume the space of $\vx$ has been shifted and transformed such that $\vx$ has mean zero and identity covariance. Consider the following symmetric and positive semidefinite matrix and its eigendecomposition,
\begin{equation} \label{eq:as_c}
\mC \;=\; \int\nabla f(\vx)\, \nabla f(\vx)^T\, \rho(\vx) \,d\vx \;=\; \mW\Lambda \mW^T \in \mathbb{R}^{m \times m},
\end{equation}
where $\Lambda=\text{diag}(\lambda_1,\dots,\lambda_m)$ and the eigenvalues are ordered from largest to smallest. Assume that for some $n < m$, $\lambda_n > \lambda_{n+1}$ (that is, the $n$th eigenvalue is \emph{strictly} greater than the $n+1$th eigenvalue). Then we can partition the eigenvectors into the first $n$, denoted $\mW_1$, and the remaining $m-n$, denoted $\mW_2$. The \emph{active subspace} is the column span of $\mW_1$, and its orthogonal complement (the column span of $\mW_2$) is the \emph{inactive subspace}. Any $\vx$ can be decomposed in terms of the active and inactive subspaces\cite{Constantine2015},
\begin{equation} \label{paramsub}
\vx \;=\; \mW_1\vy + \mW_2\vz \quad \vy \in \mathbb{R}^n, \quad \vz \in \mathbb{R}^{m-n}.
\end{equation}
We call $\vy$ the \emph{active variables} and $\vz$ the \emph{inactive variables}. Perturbing $\vy$ changes $f$ more, on average, than perturbing $\vz$; see Constantine et al.\cite{constantine2014active}. If the eigenvalues $\lambda_{n+1},\dots,\lambda_m$ are sufficiently small, then we can approximate
\begin{equation} \label{eq:gfunc}
f(\vx) \;\approx\; g(\mW_1^T\vx), 
\end{equation}
for some $g:\mathbb{R}^n\rightarrow\mathbb{R}$. In effect, this approximation enables design studies to proceed on the $n$ active variables $\vy=\mW_1^T\vx$ instead of the $m$-dimensional design variables $\vx$---i.e., dimension reduction.

With basic terms defined, we remark on the given density $\rho(\vx)$. In real engineering applications, this density function is rarely if ever known completely; thus, the engineer must choose $\rho$ such that it models the known or desired variability in the inputs $\vx$. In such cases, $\rho$ is a modeling choice, and different $\rho$'s produce different results (e.g., different active subspaces for the same $f$). The degree to which one choice of $\rho$ is better than another depends only on the suitability and appropriateness of each $\rho$ to represent the parameter variability in a particular problem. If data is available on the parameter variability, then the chosen $\rho$ should be consistent with that data. In the absence of data, the subjective choice of $\rho$ is informed by engineering intuition and design experience. There is no generic prescription in the active subspace-based toolbox for choosing $\rho$; the only constraint is that $\rho$ must be a probability density. From there, the active subspace depends on the chosen $\rho$.

The active subspace also depends on the function $f$. In general, two different functions (e.g., two different quantities of interest) contain different active subspaces, which distinguishes active subspaces from covariance-based dimension reduction techniques such as principal component analysis\cite{Jolliffe2002} (PCA) that only depend on $\rho$. Principal components can define a shape parameterization given airfoil measurements\cite{Garzon2003}, which yields the left side of the composite map \eqref{eq:composition}. Thus, it is possible to use a PCA-based shape parameterization, combined with a differentiable scalar-valued output of interest, to produce well defined active subspaces. The \emph{scalar-valued} distinction is important; other dimension reduction techniques, such as the proper orthogonal decomposition\cite{Willcox2002} (POD), seek to reduce the dimension of the state space of a large dynamical system or discretized PDE system. Although both active subspaces and POD use eigenpairs to define low-dimensional linear subspaces, these subspaces relate to different elements of a model-based design problem---active subspaces on the input parameters and POD on the state vector of the discretized flow field solutions.

\subsection{Global quadratic model for estimating active subspaces}
\label{ssec:quad}
When the gradient $\nabla f(\vx)$ is available as a subroutine (e.g., through an adjoint solver or algorithmic differentiation\cite{griewank08}), one can use Monte Carlo or other numerical integration methods to estimate $\mC$ from \eqref{eq:as_c}\cite{constantine2015computing}; the eigenpairs of the numerically estimated $\mC$ estimate those of the true $\mC$. However, many legacy solvers do not have such gradient routines. Finite differences build a local linear model at the point of the gradient evaluation and use the model's slope as the approximate gradient. But finite differences are expensive when each evaluation of $f(\vx)$ is expensive and the dimension of $\vx$ is large. We pursue a different model-based approach: (i) fit a global quadratic model with least-squares and (ii) estimate $\mC$ with the quadratic model. This approach is motivated as follows. Suppose $f(\vx)$ is well approximated by a quadratic function,
\begin{equation}
\label{eq:quadmodel}
f(\vx) \;\approx\; 
\frac{1}{2} \vx^T\mH\vx + \vx^T\vv + c,
\end{equation}
for some symmetric matrix $\mH\in\mathbb{R}^{m\times m}$, $\vv\in\mathbb{R}^m$, and $c\in\mathbb{R}$. The approximation is in the mean-squared sense with respect to the given density $\rho(\vx)$. Additionally, suppose the gradient $\nabla f(\vx)$ is well approximated as
\begin{equation}
\label{eq:quadgrad}
\nabla f(\vx) \;\approx\; \mH\vx + \vv.
\end{equation}
Again, the approximation is in the mean-squared sense. Substituting \eqref{eq:quadgrad} into \eqref{eq:as_c} and applying \eqref{eq:assump} results in the approximation
\begin{equation}
\label{eq:Capprox}
\mC 
\approx
\int\left(\mH\vx + \vv\right) \left(\mH\vx + \vv\right)^T\, \rho(\vx)\, d\vx\\
= \mH^2 + \vv\,\vv^T.
\end{equation}
The eigenpairs of $\mH^2+\vv\vv^T$ estimate those of $\mC$ in \eqref{eq:as_c},
\begin{equation}
\label{eq:Heigs}
\mH^2 + \vv\,\vv^T \;=\; \hmW\hLambda\hmW^T.
\end{equation}
The quality of the estimation depends on how well $f(\vx)$ and its gradient can be modeled by the global quadratic model \eqref{eq:quadmodel}. We emphasize that the quadratic model \eqref{eq:quadmodel} is not used as a response surface for any design or optimization problem. Its sole use is to identify important directions in the parameter space and estimate the active subspace. 

We fit the coefficients $\mH$ and $\vv$ of the quadratic model with a standard discrete least-squares approach as analyzed by Chkifa et al.\cite{Chkifa2015}. That is, we (i) draw $N$ independent samples of $\vx_i$ according to the given density $\rho$, (ii) evaluate $f_i=f(\vx_i)$ for each sample, and (iii) use least-squares with the pairs $(\vx_i, f_i)$ to fit $\mH$ and $\vv$. There are ${m+2 \choose 2}$ coefficients to fit, so we take $N$ to be 2-to-10 times ${m+2 \choose 2}$ for the least-squares fitting, which is a common rule-of-thumb. 

A related quadratic model-based approach appears in Tipireddy and Ghanem\cite{Tipireddy2014} when adapting a Hermite polynomial chaos basis for low-dimensional surrogate construction. In statistical regression, Li~\cite{Li1992} proposes using a least-squares-fit regression model to estimate the \emph{principal Hessian directions}, which are derived from eigenvectors of an averaged Hessian matrix of $f(\vx)$---in contrast to the expression in \eqref{eq:Capprox}. 

We treat the estimated eigenpairs $\hmW$, $\hLambda$ the same way we treat any other estimates of $\mC$'s eigenpairs. That is, we examine the eigenvalues in $\hLambda$. If a large gap exists between eigenvalues $\hat{\lambda}_n$ and $\hat{\lambda}_{n+1}$, then we choose $n$ to be the dimension of the active subspace. Then the first $n$ columns of $\hmW$, denoted $\hmW_1$, define the active subspace, which we can use to build a low-dimensional approximation to $f(\vx)$ as in \eqref{eq:gfunc}.

We use a bootstrap-based\cite{Efron1993} heuristic to assess the variability and uncertainty in the estimated eigenpairs from \eqref{eq:Heigs}. We apply the bootstrap to the pairs $\{(\vx_i,f_i)\}$ used to fit the quadratic model. For $k$ from 1 to $\Nboot$, 
\begin{enumerate} \label{list:bootrout}
\item let $\pi^k_j$ be an independent draw from the set $\{1,\dots,N\}$ for $j=1,\dots,N$,
\item use the pairs $\{(\vx_{\pi_j^k},f_{\pi_j^k})\}$ to fit the coefficients $\mH_k$ and $\vv_k$ of a quadratic model with least squares,
\item estimate the eigenpairs $\hmW^{(k)}$, $\hLambda^{(k)}$ of the matrix $\mH_k^2+\vv_k\vv_k^T$.
\end{enumerate}
The eigenpairs $\hmW^{(k)}k$, $\hLambda^{(k)}$ constitute one bootstrap sample. We typically choose $\Nboot$ to be between 100 and 1000, depending on the dimension of $\vx$; if $\vx$'s dimension is less than 100, then we can easily run thousands of eigenvalue decompositions from step 3 on a laptop. Note that the pairs in step 2 likely contain duplicates, which is essential in the sampling-with-replacement bootstrap method. Given the eigenvalue bootstrap samples $\{\hLambda^{(k)}\}$, we use the range over bootstrap samples to estimate uncertainty in the eigenvalues $\hLambda$. To estimate uncertainty in the $n$-dimensional active subspace, let $\hmW_{1}^{(k)}$ contain the first $n$ columns of the bootstrap sample of eigenvectors $\hmW^{(k)}$, and compute the subspace error\cite{Golub1996},
\begin{equation}
\label{eq:booterror}
\hat{e}_n^{(k)} \;=\;
\left\| \left(\hmW_{1}^{(k)}\right)\left(\hmW_{1}^{(k)}\right)^T
- \left(\hmW_{1}\right)\left(\hmW_{1}\right)^T \right\|_2,
\end{equation}
where the norm is the matrix 2-norm. The mean and standard error of the set $\{e_n^{(k)}\}$ estimates the error and uncertainty in the active subspace defined by $\hmW_1$ from \eqref{eq:Heigs}. 

There are two important caveats to the bootstrap-based heuristic for assessing uncertainty in the computed quadratic model-based active subspace. First, the pairs $\{(\vx_i,f_i)\}$ (i.e., the data) used to fit the quadratic model do not satisfy the assumptions needed to allow a proper statistical interpretation of the bootstrap-based error estimates. Proper interpretations require structured noise on the $f_i$'s; when $f_i$ is the output of a deterministic computer simulation, such noise models are not valid. Second, the bootstrap-based error estimates only quantify the error due to finite sampling when estimating the quadratic model. They do not include errors from the motivating quadratic model assumption \eqref{eq:quadmodel} and \eqref{eq:quadgrad}. Despite these caveats, the bootstrap-based error estimates give some indication of variability in the eigenpairs estimated in \eqref{eq:Heigs}. 

In section \ref{sec:results}, we show numerical evidence of convergence as $N$ increases for active subspaces derived from the composite map \eqref{eq:composition} with both PARSEC and CST as the underlying airfoil shape parameterizations. 

\section{Active subspaces in PARSEC and CST parameterizations}
\label{sec:results}

The goal of the following numerical experiments is to study how the active subspaces of lift and drag change as the shape parameterization changes. That is, if we change the first component of the composite map \eqref{eq:composition}, $\vx\rightarrow s_U,s_L$, (i) how do the active subspaces change for the composite map and (ii) how do insights into the design problem change? In particular, we compare active subspaces derived using both PARSEC and CST parameterizations for the airfoil surfaces. 

In each case, we must set the density function $\rho(\vx)$ needed to define the active subspaces. We want to ensure that all parameters $\vx$ in the support of $\rho$ produce shapes that satisfy the necessary conditions from section \ref{sec:airfoil}; our approach is somewhat conservative. For each parameterization, we first find the parameters that produce a shape nearest the NACA 0012 using nonlinear least-squares, which defines the respective center of the parameter space. We then set bounds on the parameters to be $\pm20$\% of the center. The density $\rho$ is uniform on the hyperrectangle defined by these ranges. 


Given a shape, we estimate the quantities of interest lift and drag at fixed $0^{\circ}$ angle of attack and NIST conditions (ambient temperature $0^{\circ}$ C, ambient pressure $101.325$ KPa, Mach $0.8$) using the Stanford University Unstructured (SU2) CFD solver\cite{economon2016}, which solves Euler equations for compressible flow on an unstructured mesh. We executed each run on the Colorado School of Mines Mio cluster. The data supporting these results can be found on GitHub (\url{https://github.com/zgrey/ASAP.git}).

\subsection{PARSEC shape parameterization results}
The bounds for the PARSEC parameters (see Table \ref{tab:PARSECparms}) are in Table \ref{PARStable}. These bounds are $\pm20$\% of fitted values for a NACA 0012 airfoil. To generate the pairs needed for the quadratic model-based approximation of the active subspaces, we draw $N=6500$ independent samples from the uniform density on the hyperrectangle defined by the ranges in Table \ref{PARStable}. Each sample produces a shape that leads to a CFD solve and lift and drag coefficients, $C_l$ and $C_d$, respectively. The pairs $\{(\vx_i, C_l(\vx_i))\}$ are used to estimate the active subspaces for lift as a function of PARSEC parameters following the procedure outlined in section \ref{ssec:quad}. And the pairs $\{(\vx_i, C_d(\vx_i))\}$ are similarly used to estimate the active subspaces for drag.

\begin{table}[H]
\caption{Bounds on PARSEC parameters (see Table \ref{tab:PARSECparms}) chosen by first finding parameters that produce the NACA 0012 and then setting ranges to be $\pm20$\% of the optimized parameters.}
\begin{center}
\begin{tabular}{l|l|l}
Parameter & Lower & Upper\\
\hline
	$x_1$ & 0.242 & 0.363\\
	$x_2$ & 0.242 & 0.363\\
	$x_3$ & 0.048 & 0.072\\
	$x_4$ & -0.072 & -0.048\\
	$x_5$ & -0.004 & 0.004\\
	$x_6$ & 0.008 & 0.012\\
	$x_7$ & -2.223 & -3.335\\
	$x_8$ & 7.40 & 11.10\\
	$x_9$ & -0.400 & -0.600\\
	$x_{10}$ & 0.400 & 0.600\\
	$x_{11}$ & 0.012 & 0.018\\
\end{tabular}
\end{center}
\label{PARStable}
\end{table}

\begin{figure}[H] 
\centering
\subfloat[Lift]{
\includegraphics[width=0.45\textwidth]{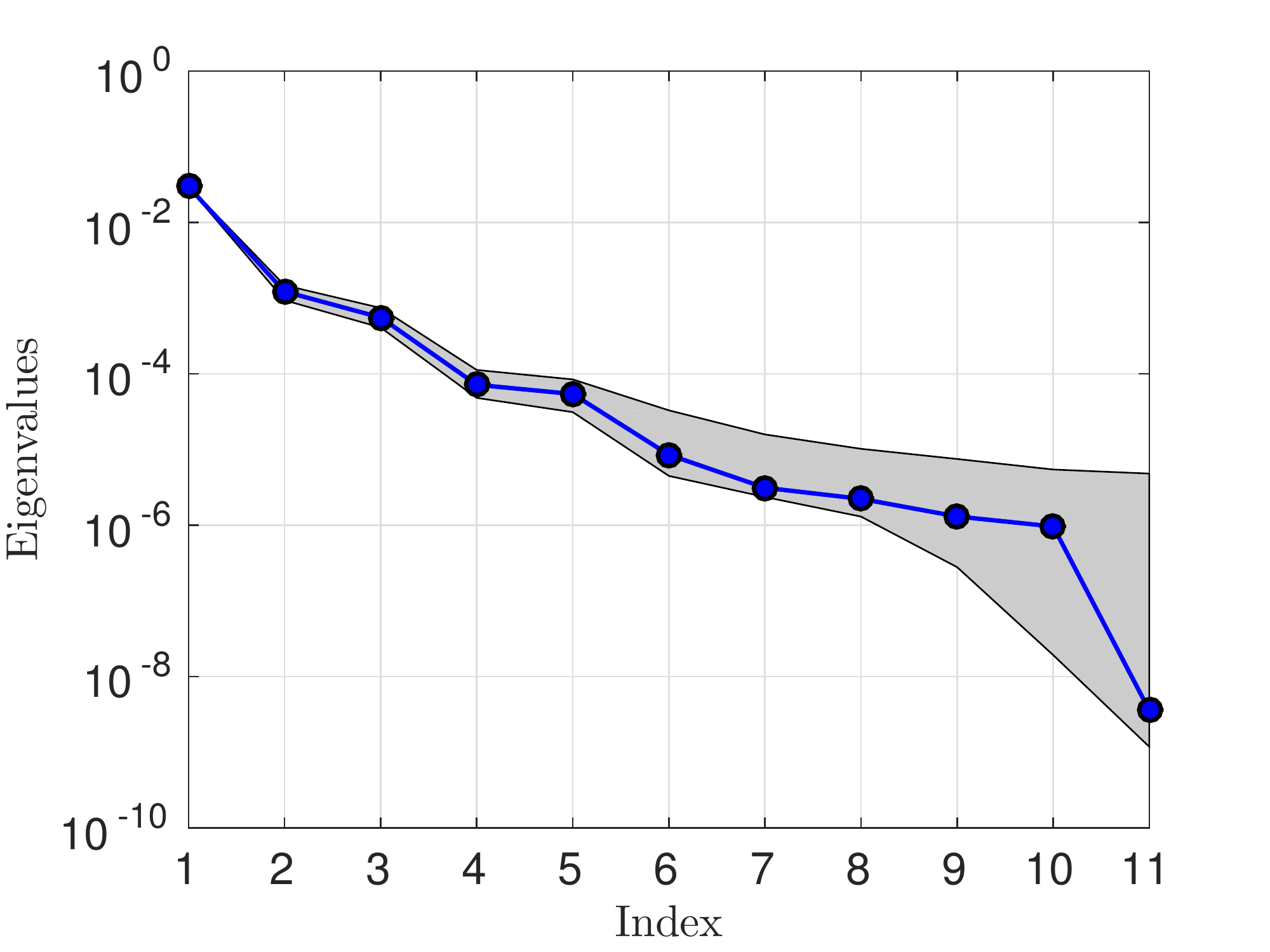}
}
\hfil
\subfloat[Drag]{
\includegraphics[width=0.45\textwidth]{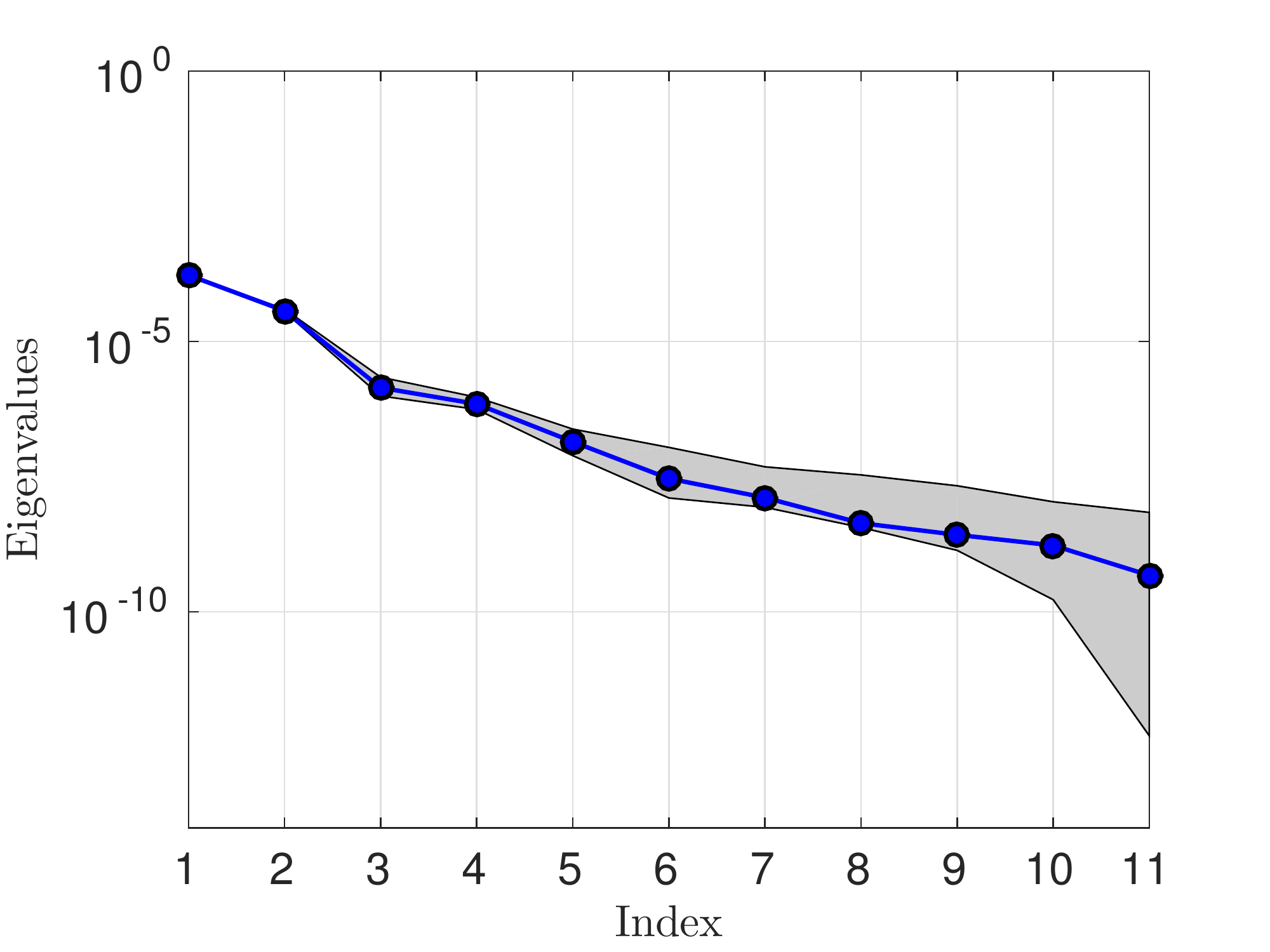}
}
\caption{Eigenvalues from \eqref{eq:Heigs} for lift (left) and drag (right) functions of the PARSEC parameters from Table \ref{PARStable}. }
\label{PARSeig}
\end{figure}

\begin{figure}[H]
\centering
\subfloat[Lift]{
\includegraphics[width=0.45\textwidth]{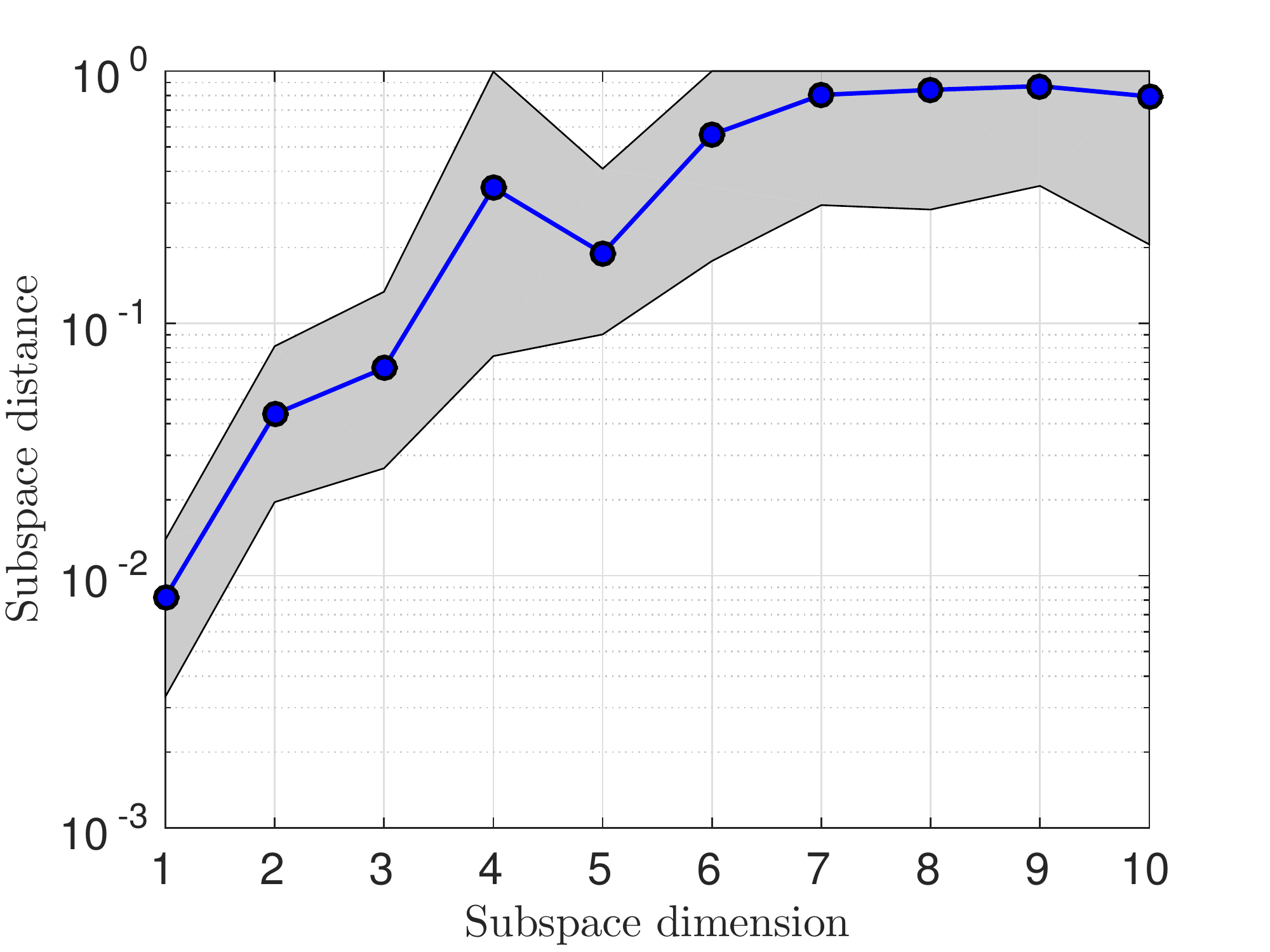}
}
\hfil
\subfloat[Drag]{
\includegraphics[width=0.45\textwidth]{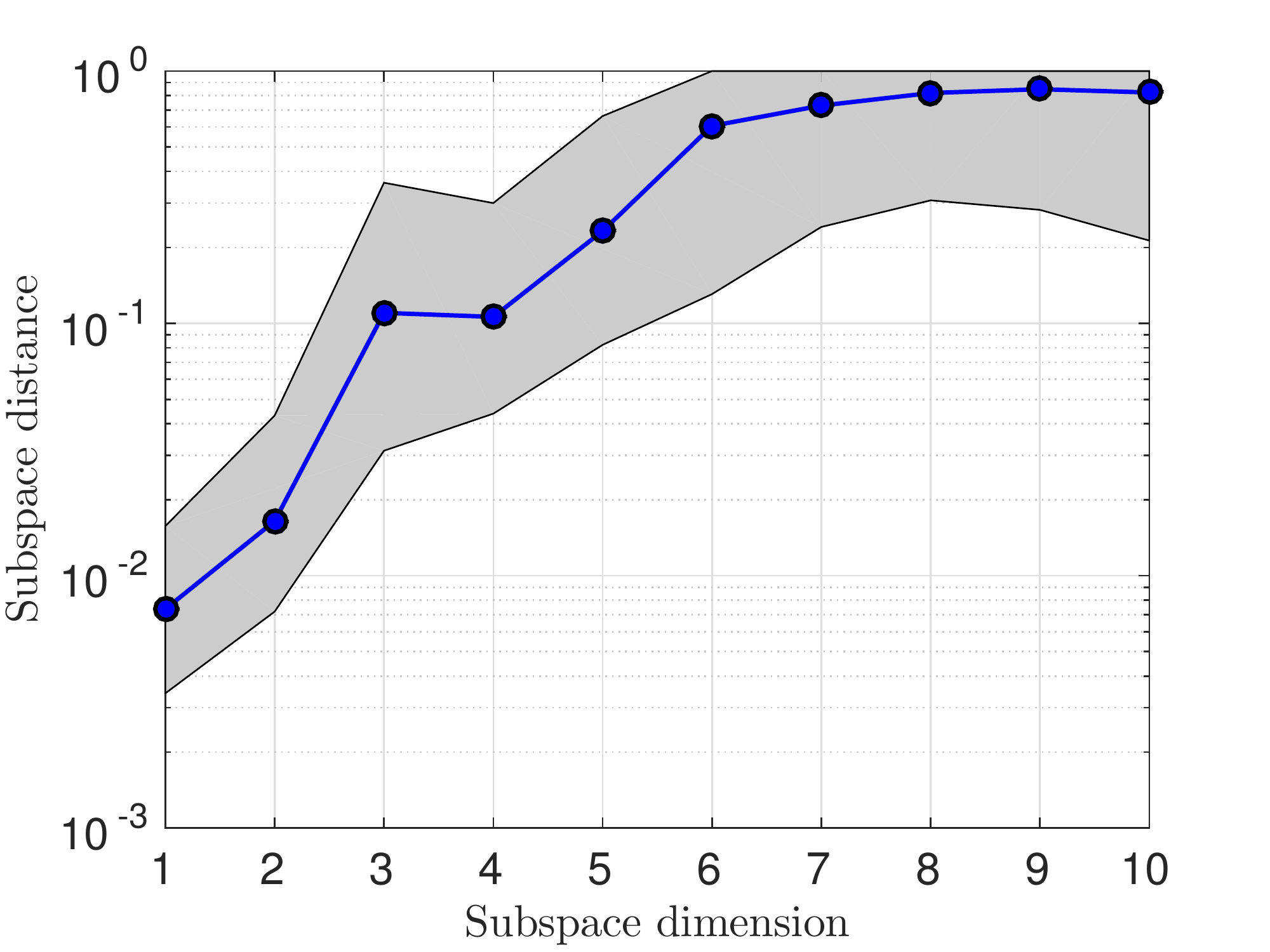}
}
\caption{Subspace error estimates (means and ranges) from \eqref{eq:booterror} for active subspaces estimated with the quadratic model-based approach (section \ref{ssec:quad}) for lift (left) and drag (right) using the PARSEC shape parameterization.}
\label{PARSerr}
\end{figure}

Figure \ref{PARSeig} shows the 11 eigenvalues $\hLambda$ from \eqref{eq:Heigs} for lift and drag on a logarithmic scale along with the bootstrap ranges that estimate uncertainty (see section \ref{ssec:quad}); we use $\Nboot=500$ bootstrap samples. For lift, there is a large gap between the first and second eigenvalues, which is evidence of an exploitable one-dimensional active subspace. For drag, a relatively large gap occurs between the second and third eigenvalues, which suggests an exploitable two-dimensional active subspace. Figure \ref{PARSerr} shows, for both lift and drag, the means and ranges for the bootstrap-based subspace error estimates \eqref{eq:booterror} as a function of the active subspace dimension. Notice that the estimated error increases as the subspace dimension increases. This is related to the decreasing gaps between subsequent eigenvalues from Figure \ref{PARSeig}. As the gap between two eigenvalues decreases, estimating the eigenspace whose dimension is associated with the smaller eigenvalue becomes more difficult for numerical procedures\cite{Stewart1973}. In the extreme case, if $\lambda_n=\lambda_{n+1}$, then there is no eigenspace of dimension $n$. 

\begin{figure}[H]
\centering
\subfloat[Lift]{\label{fig:parsec-lift-evecs}
\includegraphics[width=0.45\textwidth]{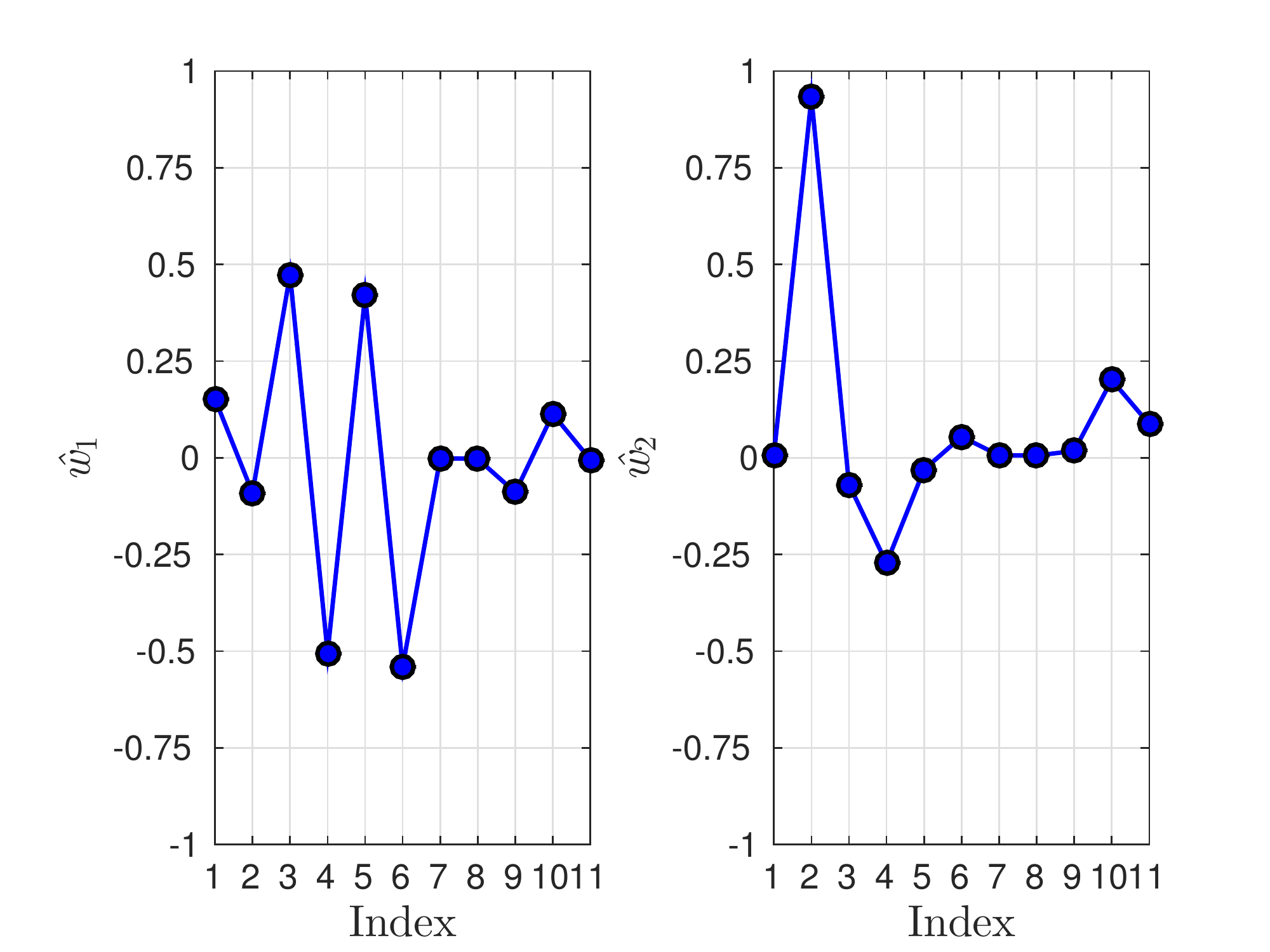}
}
\hfil
\subfloat[Drag]{\label{fig:parsec-drag-evecs}
\includegraphics[width=0.45\textwidth]{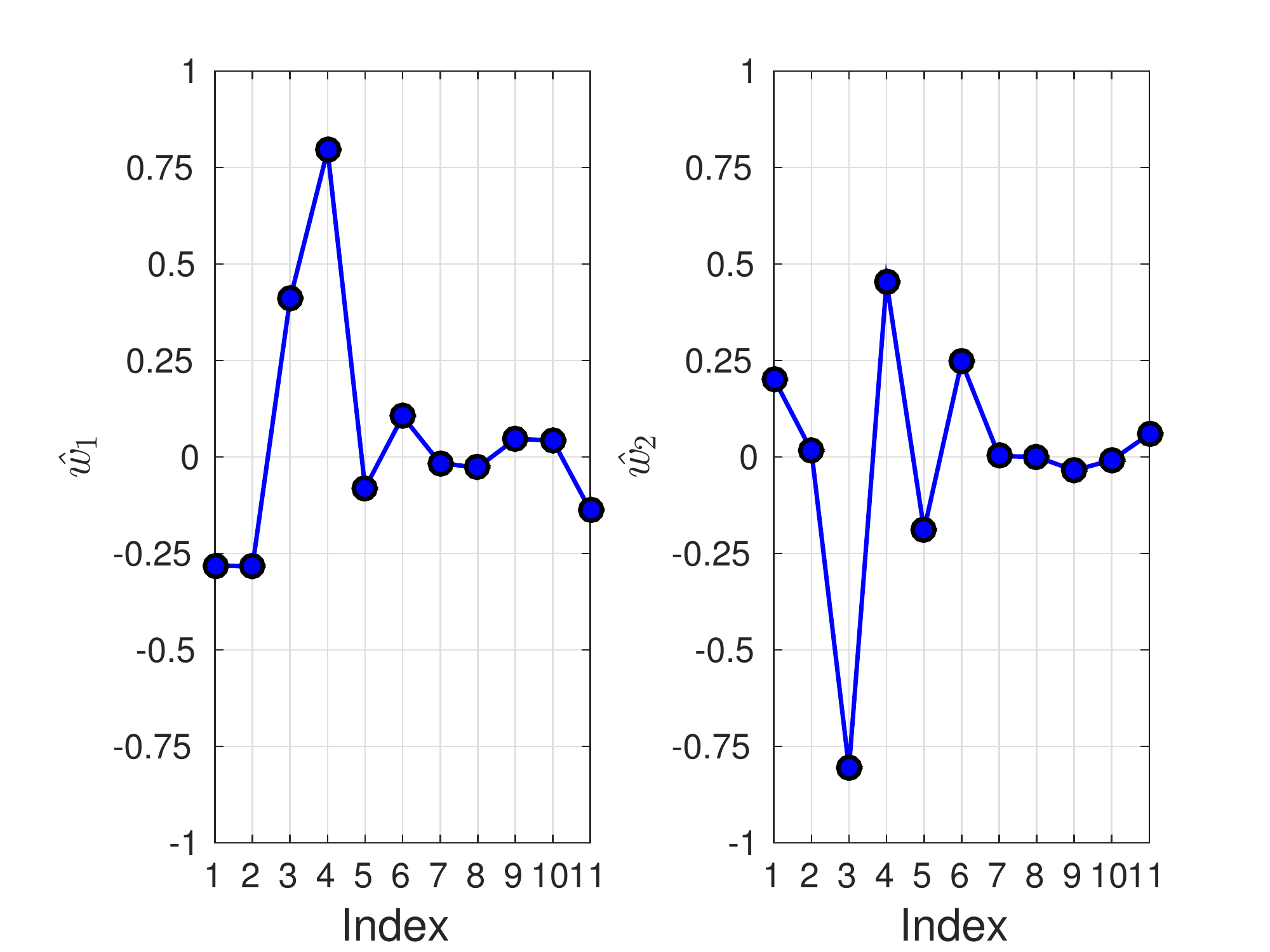}
}
\caption{Components of the first two eigenvectors from \eqref{eq:Heigs} for lift (left) and drag (right) using the PARSEC shape parameterization.}
\label{PARSvec}
\end{figure}

Figure \ref{PARSvec} shows components of the eigenvectors $\hmW$ from \eqref{eq:Heigs}. Figure \ref{fig:parsec-lift-evecs} shows the first two eigenvectors associated with lift. And Figure \ref{fig:parsec-drag-evecs} shows the first two eigenvectors associated with drag. These eigenvector components can be used to compute sensitivity metrics for the PARSEC parameters with respect to the outputs of interest\cite{Constantine2017}. The indices on the horizontal axis labels correspond to the labels in Table \ref{tab:PARSECparms}. A large eigenvector component implies that output is more sensitive to changes in the corresponding parameter. The eigenvector components suggest that lift is most influenced on average by $x_3$, $x_4$, $x_5$, and $x_6$, which affect the upper and lower surface coordinates and trailing edge deflections. Additionally, the eigenvector components for upper and lower surface coordinates are nearly equal and opposite ($\pm0.5$). This suggests that increasing the upper surface coordinate should be combined with reducing the lower surface coordinate to change lift the most. These coordinate changes lead to a change in asymmetry of the shape (camber). The eigenvector components for drag (Figure \ref{fig:parsec-drag-evecs}) are largest for $x_3$ and $x_4$, corresponding to the upper and lower surface coordinates, respectively. These large components have the same sign, which suggests that simultaneous changes in the upper and lower surface coordinate affect drag more, on average. Simultaneously perturbing $x_3$ and $x_4$ changes thickness. This agrees with the design intuition that changes in thickness should influence drag the most, on average. Perturbing $\vx$ along the second eigenvector for drag changes camber.

Given (i) the eigenvectors $\hmW$ from \eqref{eq:Heigs} and (ii) the lift data $\{(\vx_i, C_l(\vx_i))\}$ and drag data $\{(\vx_i,C_d(\vx_i))\}$, we can create useful visualizations called \emph{shadow plots} that elucidate the relationship between the shape parameters $\vx$ and the quantities of interest $C_l$ and $C_d$. These shadow plots are motivated by \emph{sufficient summary plots} developed by Cook\cite{Cook1998} for sufficient dimension reduction. Given a data set and basis vectors defining a subspace of the inputs, the shadow plot's construction is identical to the sufficient summary plot---but their interpretations differ. Proper interpretation of the sufficient summary plot relies on the deep theory of sufficient dimension reduction in statistical regression that justifies the label \emph{sufficient} for the plot. In contrast, the shadow plot is a qualitative tool with no notion of statistical sufficiency. Interpreting the shadow plot requires the modeler's subjective judgment. See Chapter 1 in Constantine\cite{Constantine2015} for more information on these plots using data derived from computational science simulation models. 

Let $\hvw_1$ be the first column of $\hmW$. The one-dimensional shadow plot plots the output (e.g., $C_l(\vx_i)$) on the vertical axis versus a linear combination of the associated input vector $\hvw_1^T\vx_i$ on the horizontal axis, where the weights of the linear combination are the components of the vector $\hvw_1$. If a one-dimensional summary plot shows a univariate relationship between the linear combination of inputs $\hvw_1^T\vx$ and the output, then an approximation of the form \eqref{eq:gfunc} with one vector is appropriate. If a univariate relationship is not apparent, then one can make a two-dimensional shadow plot. Let $\hvw_2$ be the second column of $\hmW$. The two-dimensional shadow plot puts the outputs on the $z$-axis, the linear combination of inputs $\hvw_1^T\vx_i$ on the $x$-axis, and the linear combination of inputs $\hvw_2^T\vx_i$ on the $y$-axis. If the output appears to be a function of these two linear combinations, then we can approximate the output as in \eqref{eq:gfunc} with two vectors. If the two-dimensional summary plot does not reveal a useful bivariate relationship, then the approximation \eqref{eq:gfunc} may be appropriate with more than two vectors, but metrics other than the shadow plots are needed to assess the approximation quality (e.g., residuals). 

\begin{figure}[H]
\centering
\subfloat[Lift]{
\includegraphics[width=0.45\textwidth]{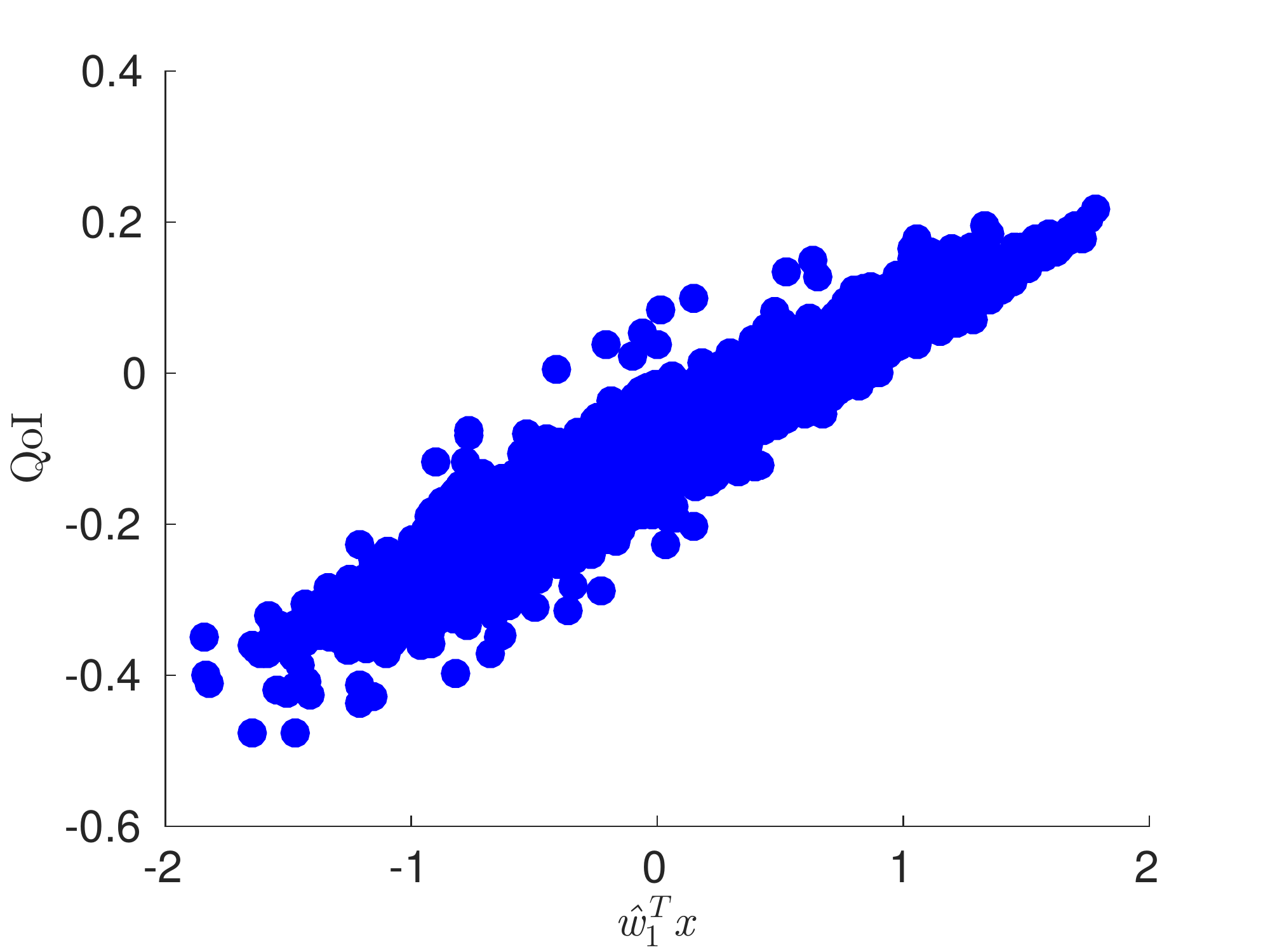}
}
\hfil
\subfloat[Drag]{
\includegraphics[width=0.45\textwidth]{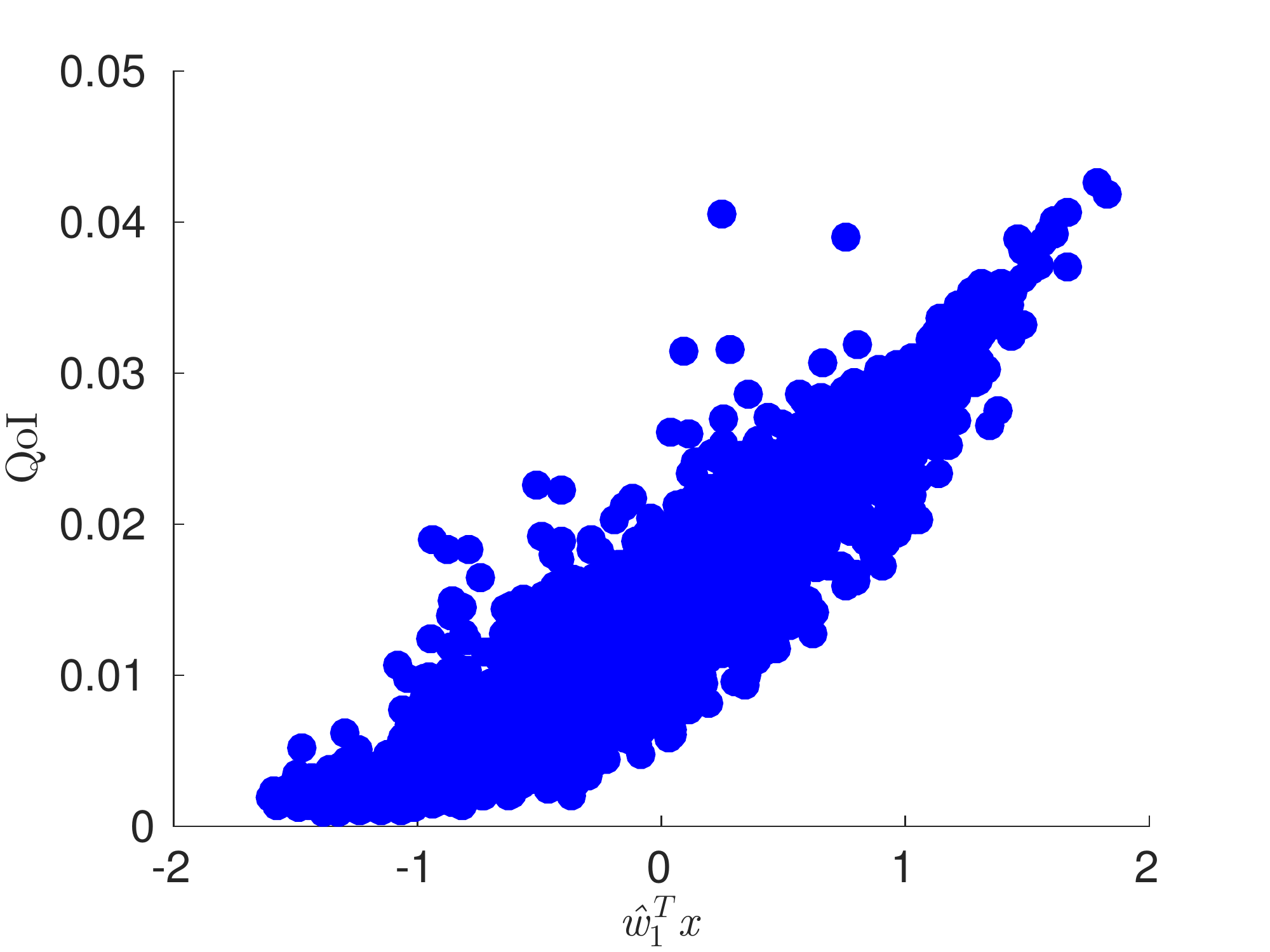}
}
\caption{One-dimensional shadow plots of lift (left) and drag (right) using the PARSEC parameterization with the first eigenvector from \eqref{eq:Heigs}.}
\label{1DPARS}
\end{figure}

\begin{figure}[H]
\centering
\subfloat[Lift]{
\includegraphics[width=0.45\textwidth]{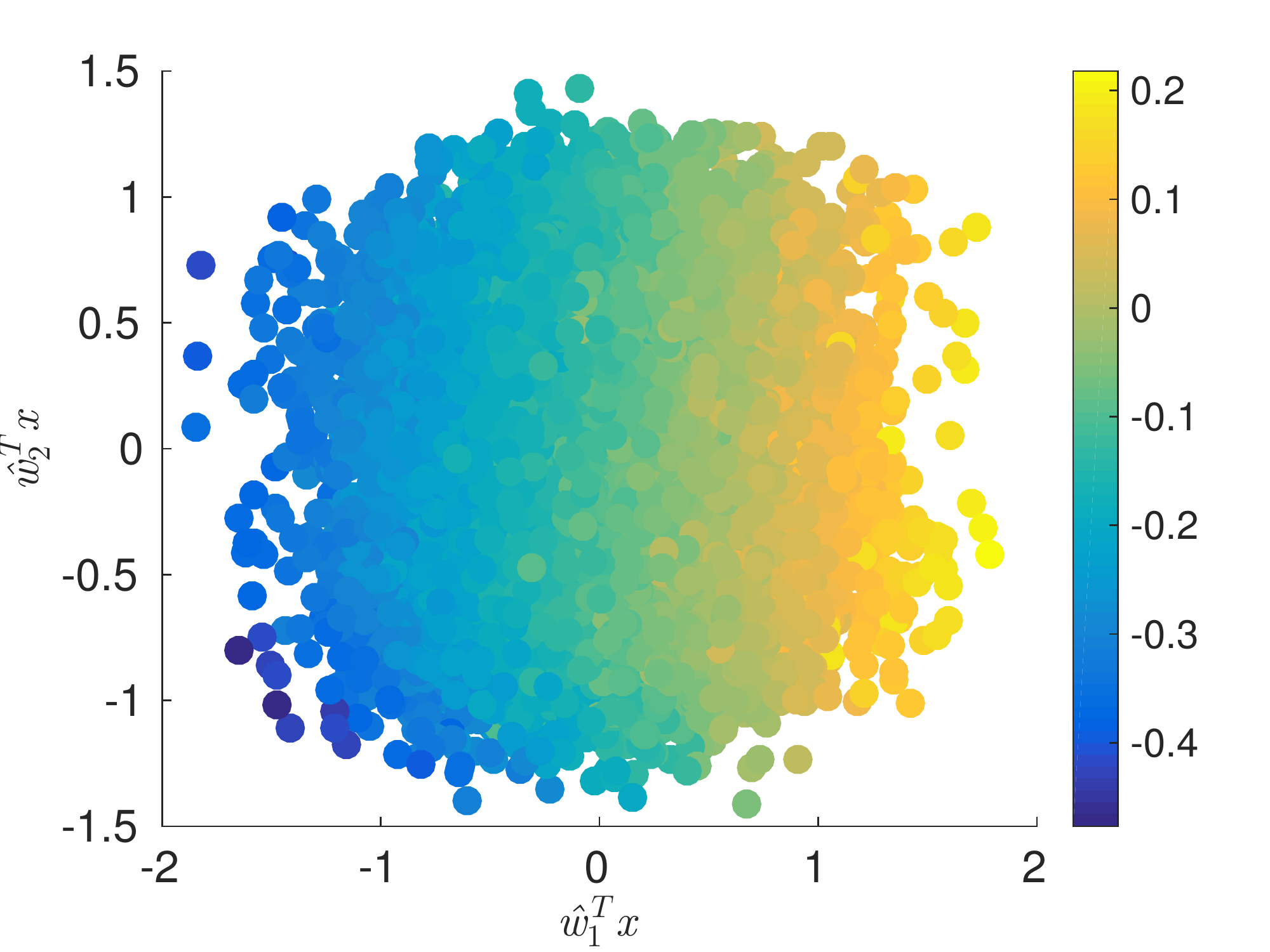}
}
\hfil
\subfloat[Drag]{
\includegraphics[width=0.45\textwidth]{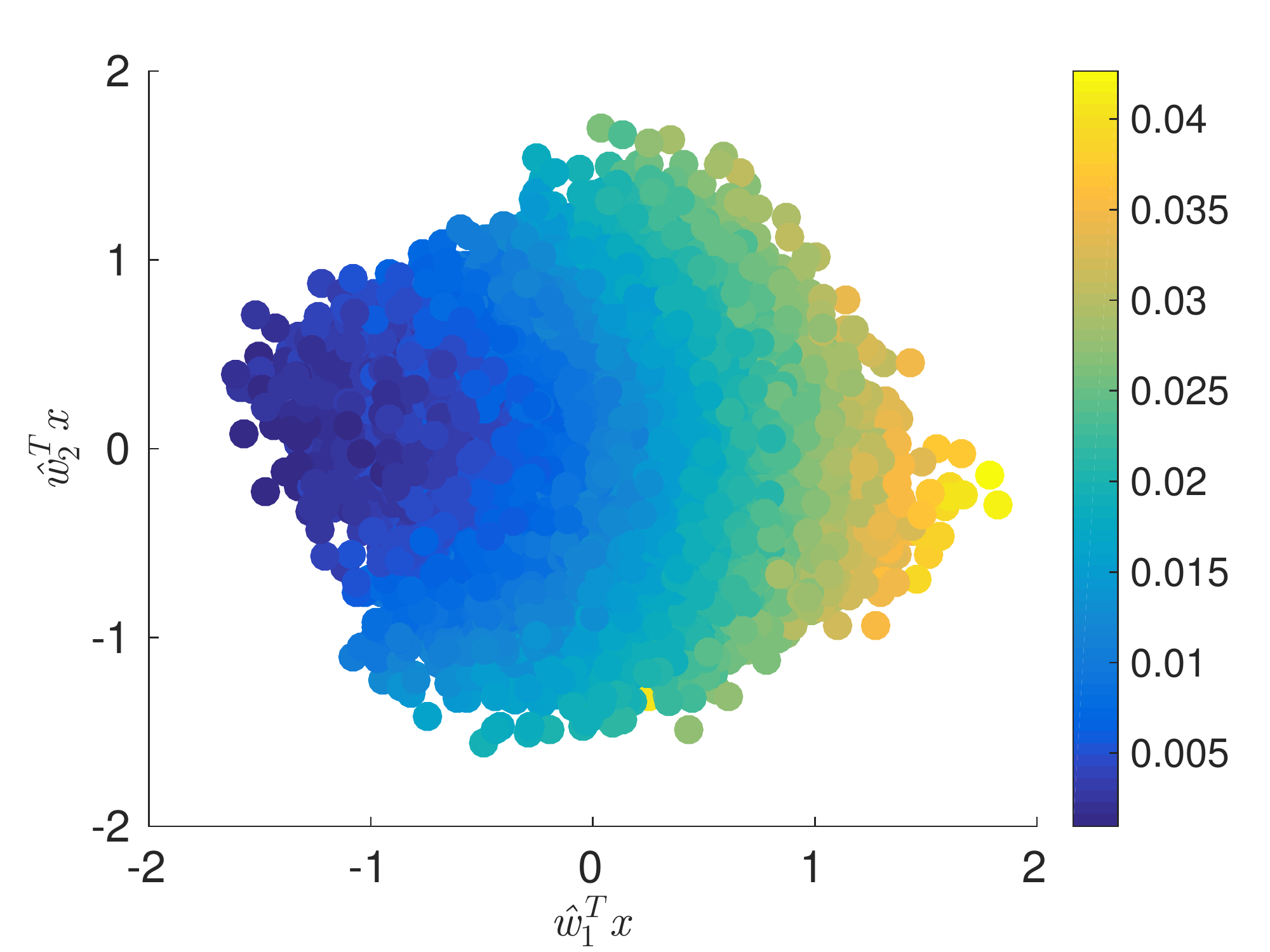}
}
\caption{Two-dimensional shadow plots of lift (left) and drag (right) using the PARSEC parameterization with the first two eigenvectors from \eqref{eq:Heigs}.}
\label{2DPARS}
\end{figure}

Figure \ref{1DPARS} shows the one-dimensional shadow plots for lift and drag using the PARSEC parameterization with all 6500 samples used to estimate the quadratic model-based active subspaces. Note the strong monotonic relationship in both cases. Such a relationship may be useful if one wishes to estimate the range of outputs over the range of inputs; see Constantine et al.\cite{Constantine2015exploiting} for an example of this in a hypersonic scramjet model. However, the relationship is far from univariate. Figure \ref{2DPARS} shows the two-dimensional shadow plots for the lift and drag. The bivariate relationship for drag is much more apparent. 

\subsection{CST shape parameterization results}
We follow the same procedure to set $\rho$ for the CST parameterization. First find values for the CST parameters such that the resulting shape matches a NACA 0012 airfoil. Then choose bounds on the parameters to be $\pm20$\% of the fitted values, and set $\rho$ to be a uniform density on the resulting hyper-rectangle. Table \ref{CSTtable} shows the bounds for the CST parameters. We draw $N=6500$ samples of the $m=10$ CST parameters independently according to $\rho$. For each parameter sample, we construct the airfoil, solve for the flow fields with SU2\cite{economon2016}, and compute the outputs of interest (lift and drag). The resulting pairs are used to estimate the quadratic model-based active subspaces for lift and drag as described in section \ref{ssec:quad}. 

\begin{table}[H]
\caption{Bounds on CST parameters (see section \ref{sec:CST}) chosen by first finding parameters that produce the NACA 0012 and then setting ranges to be $\pm20$\% of the optimized parameters.}
\begin{center}
\begin{tabular}{l|l|l}
Parameter & Lower & Upper\\
\hline
	$x_1$ & 0.12 & 0.18\\
	$x_2$  & 0.8 & 1.2\\
	$x_3$ & 0.8 & 1.2\\
	$x_4$ & 0.8 & 1.2\\
	$x_5$ & 0.8 & 1.2\\
	$x_6$ & -0.12 & -0.18\\
	$x_7$ & 0.8 & 1.2\\
	$x_8$ & 0.8 & 1.2\\
	$x_9$ & 0.8 & 1.2\\
	$x_{10}$ & 0.8 & 1.2\\
\end{tabular}
\end{center}
\label{CSTtable}
\end{table}

\begin{figure}[H]
\centering
\subfloat[Lift]{
\includegraphics[width=0.45\textwidth]{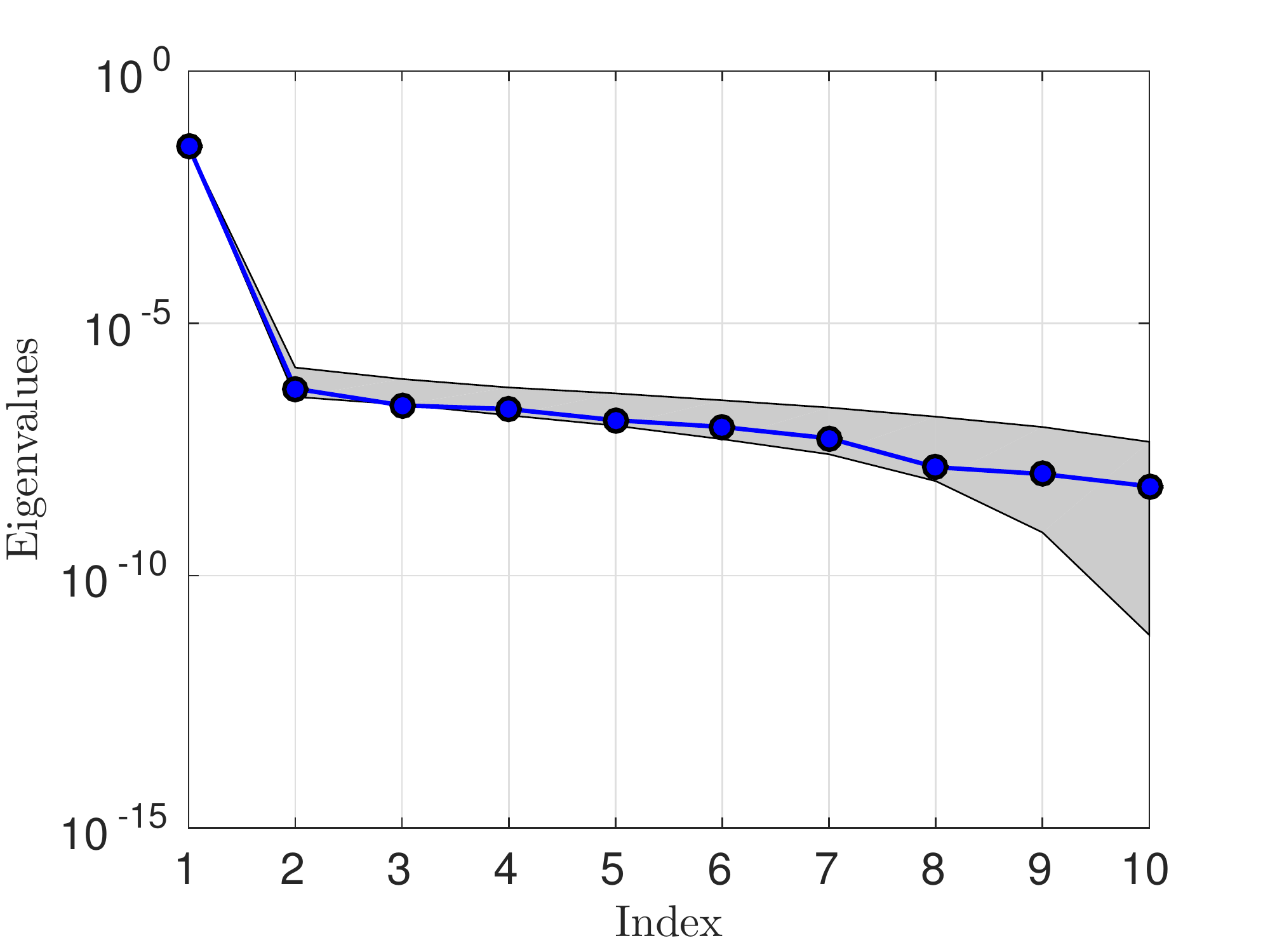}
}
\hfil
\subfloat[Drag]{
\includegraphics[width=0.45\textwidth]{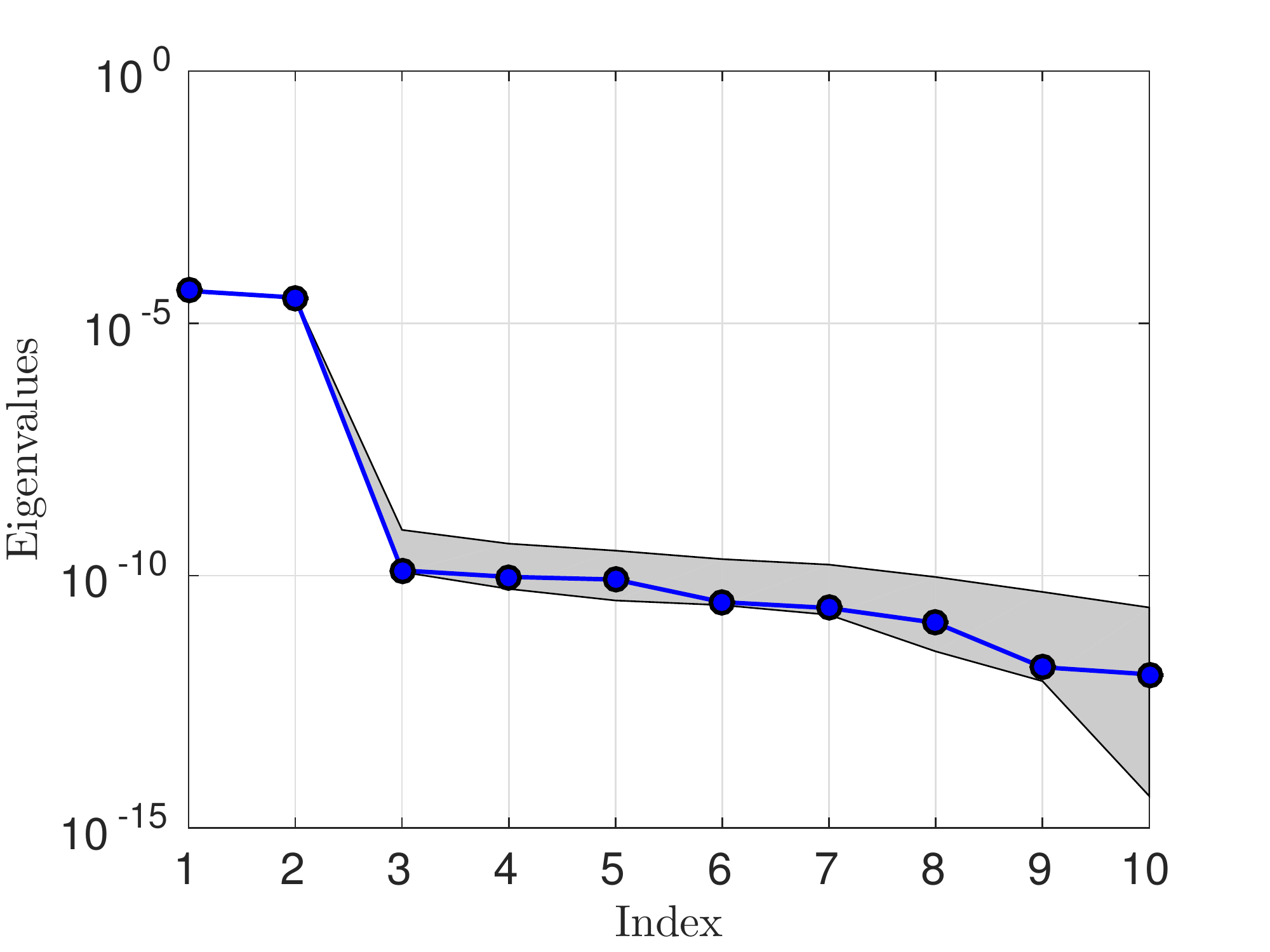}
}
\caption{Eigenvalues from \eqref{eq:Heigs} for lift (left) and drag (right) functions of the CST parameters from Table \ref{CSTtable}.}
\label{CSTeig}
\end{figure}

\begin{figure}[H]
\centering
\subfloat[Lift]{\label{fig:cst-lift-evecs}
\includegraphics[width=0.45\textwidth]{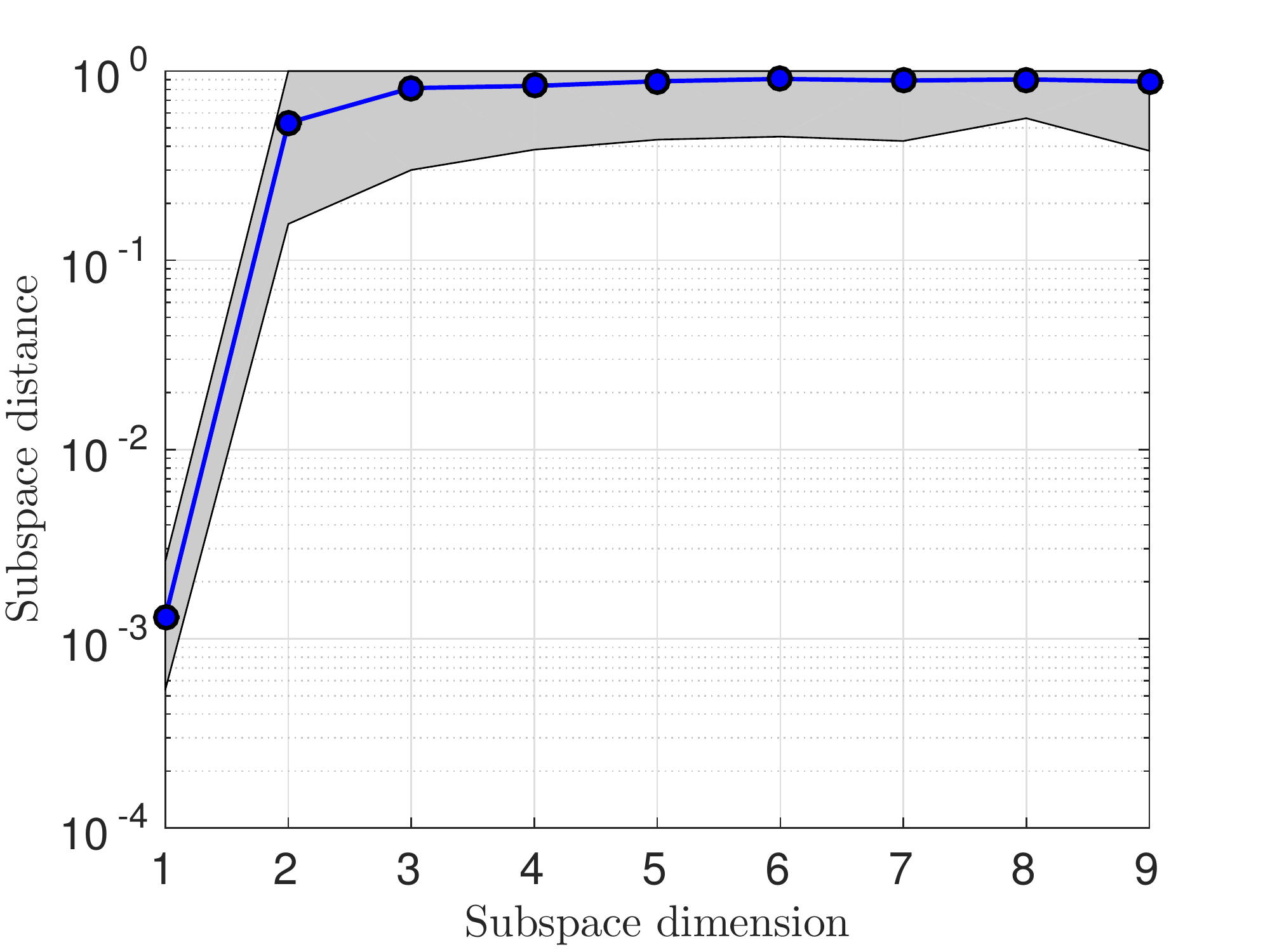}
}
\hfil
\subfloat[Drag]{\label{fig:cst-drag-evecs}
\includegraphics[width=0.45\textwidth]{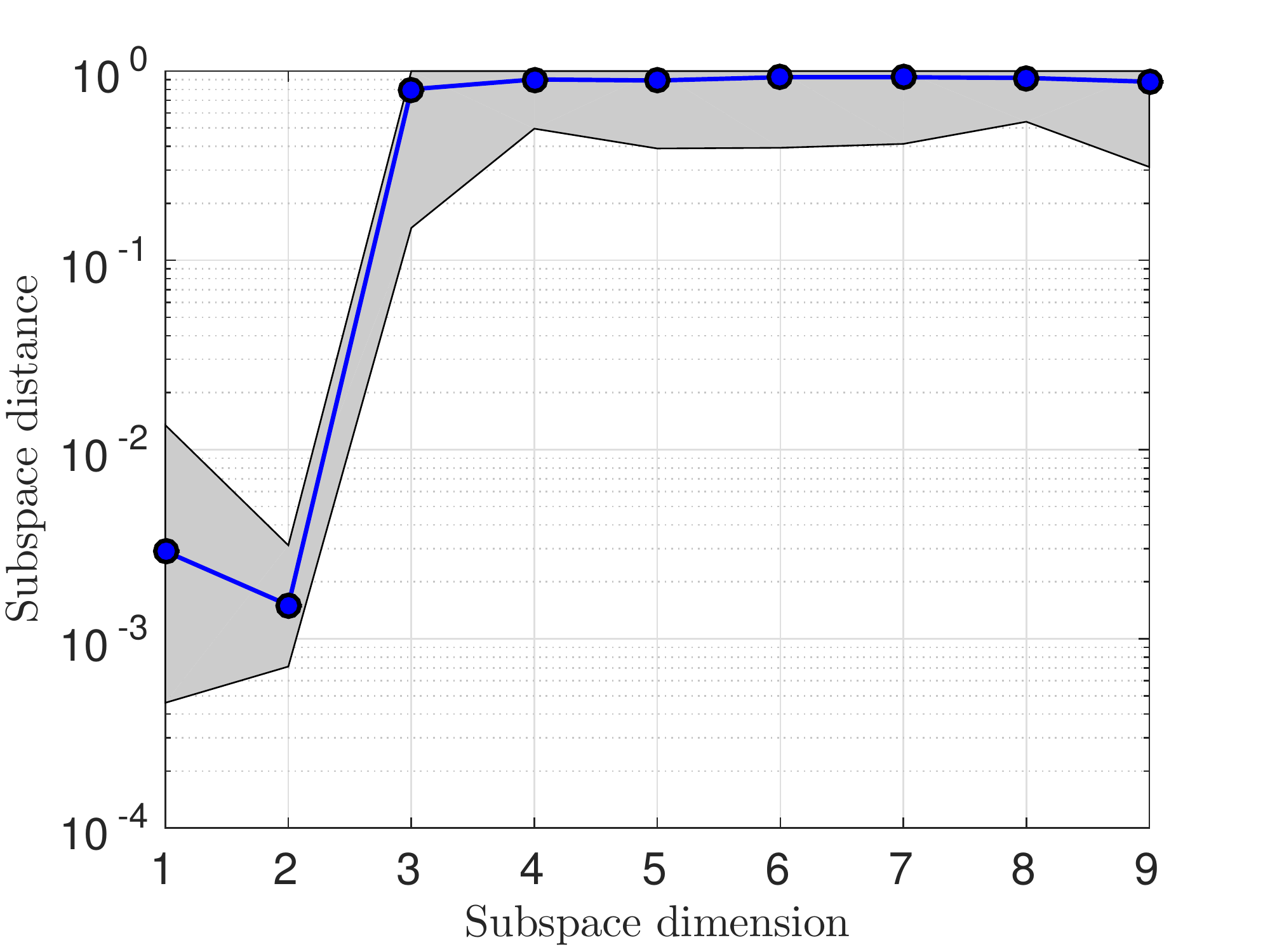}
}
\caption{Subspace error estimates (means and ranges) from \eqref{eq:booterror} for active subspaces estimated with the quadratic model-based approach (section \ref{ssec:quad}) for lift (left) and drag (right) using CST shape parameteization.}
\label{CSTerr}
\end{figure}

Figure \ref{CSTeig} shows the eigenvalues from \eqref{eq:Heigs} with bootstrap ranges for lift and drag. In contrast to the eigenvalues using the PARSEC parameterization from Figure \ref{PARSeig}, the eigenvalues using CST show very pronounced gaps (i) between the first and second eigenvalues for lift and (ii) between the second and third eigenvalues for drag. These gaps suggest that a one-dimensional active subspace is useful for approximating lift, and a two-dimensional active subspace is useful for approximating drag. These conclusions are consistent with the conclusions we drew from the PARSEC-parameterized results. Figure \ref{CSTerr} shows the bootstrap-based estimated subspace errors due to finite sampling as the dimension of the subspace increases. Similar to Figure \ref{PARSerr}, the subspace errors are inversely proportional to the associated eigenvalue gaps in Figure \ref{PARSeig}, which is consistent with well-known theory for numerical approximation of eigenspaces\cite{Stewart1973}. Thus, the errors in the estimated one-dimensional active subspace for lift and the one- and two-dimensional active subspaces for drag are small enough to be useful for design optimization. 

\begin{figure}[H]
\centering
\subfloat[Lift]{
\includegraphics[width=0.45\textwidth]{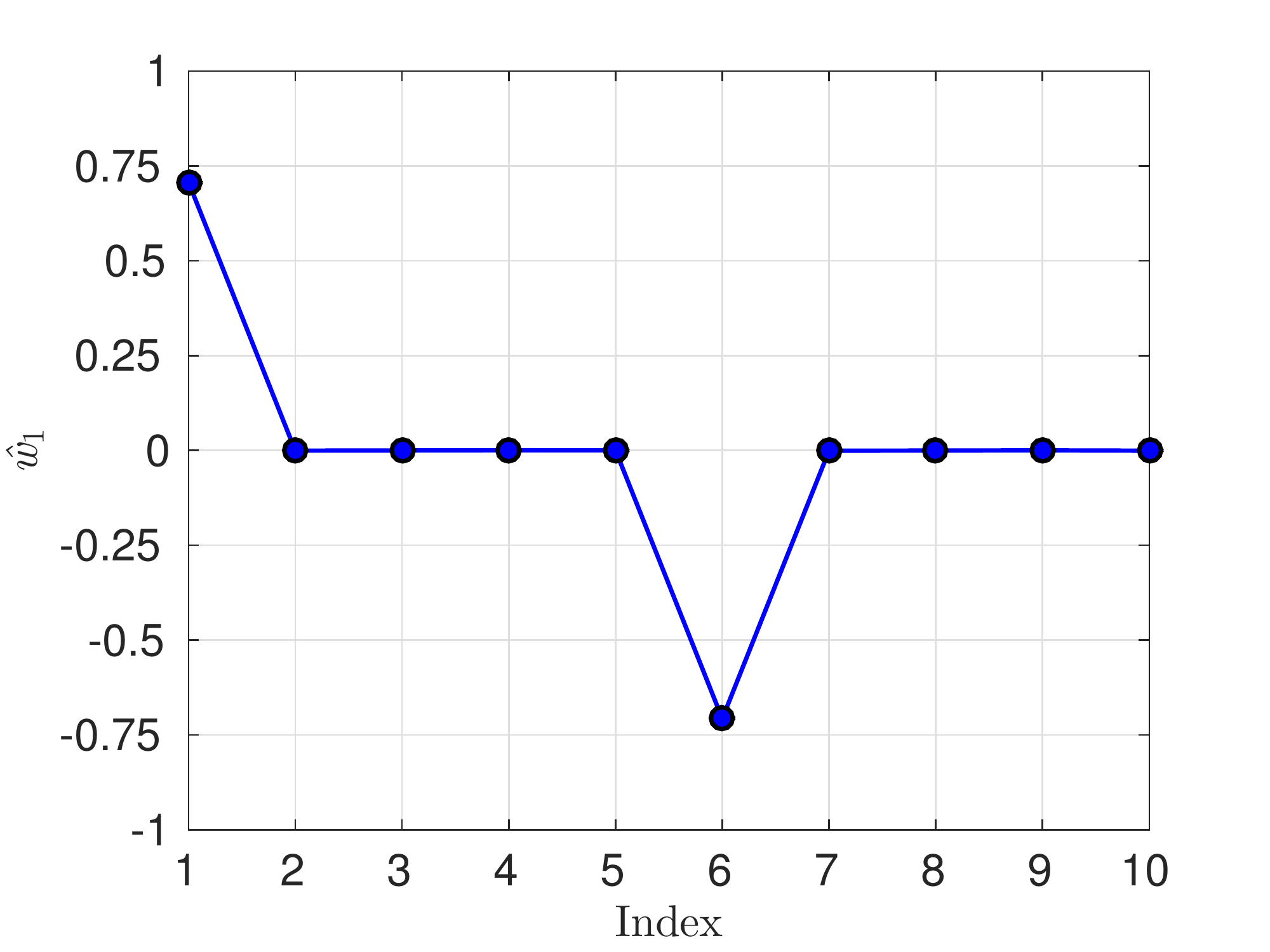}
}
\hfil
\subfloat[Drag]{
\includegraphics[width=0.45\textwidth]{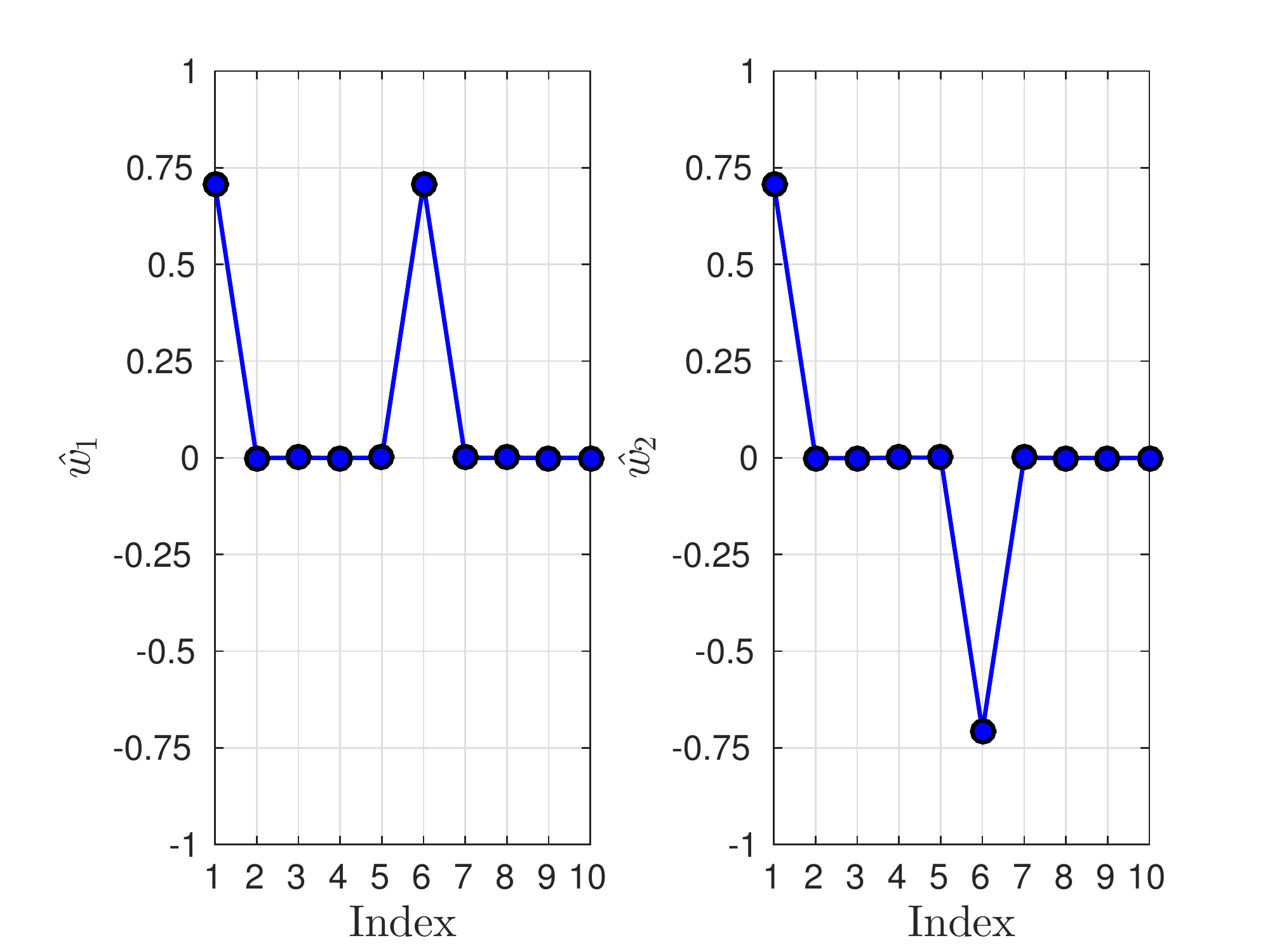}
}
\caption{Components of the first two eigenvectors from \eqref{eq:Heigs} for lift (left) and drag (right) using the CST shape parameterization.}
\label{CSTvec}
\end{figure}

Figure \ref{CSTvec} shows (i) the first eigenvector from \eqref{eq:Heigs} for lift and (ii) the first two eigenvectors for drag using the CST shape parameterization. The nonzero eigenvector components correspond to the leading coefficient in the CST series \eqref{eq:cstseriest} for each of the upper and lower surfaces---index 1 for the upper surface and index 6 for the lower surface. The eigenvector components suggest that the active subspaces are defined purely by these leading coefficients---both for lift and drag. In other words, CST parameters $x_2,\dots,x_5,x_7,\dots, x_{10}$ have no impact on either lift or drag (for the chosen $\rho(\vx)$). Additionally, the eigenvector for lift is identical to the second eigenvector for drag.

\begin{figure}[H]
\centering
\subfloat[Lift]{
\includegraphics[width=0.45\textwidth]{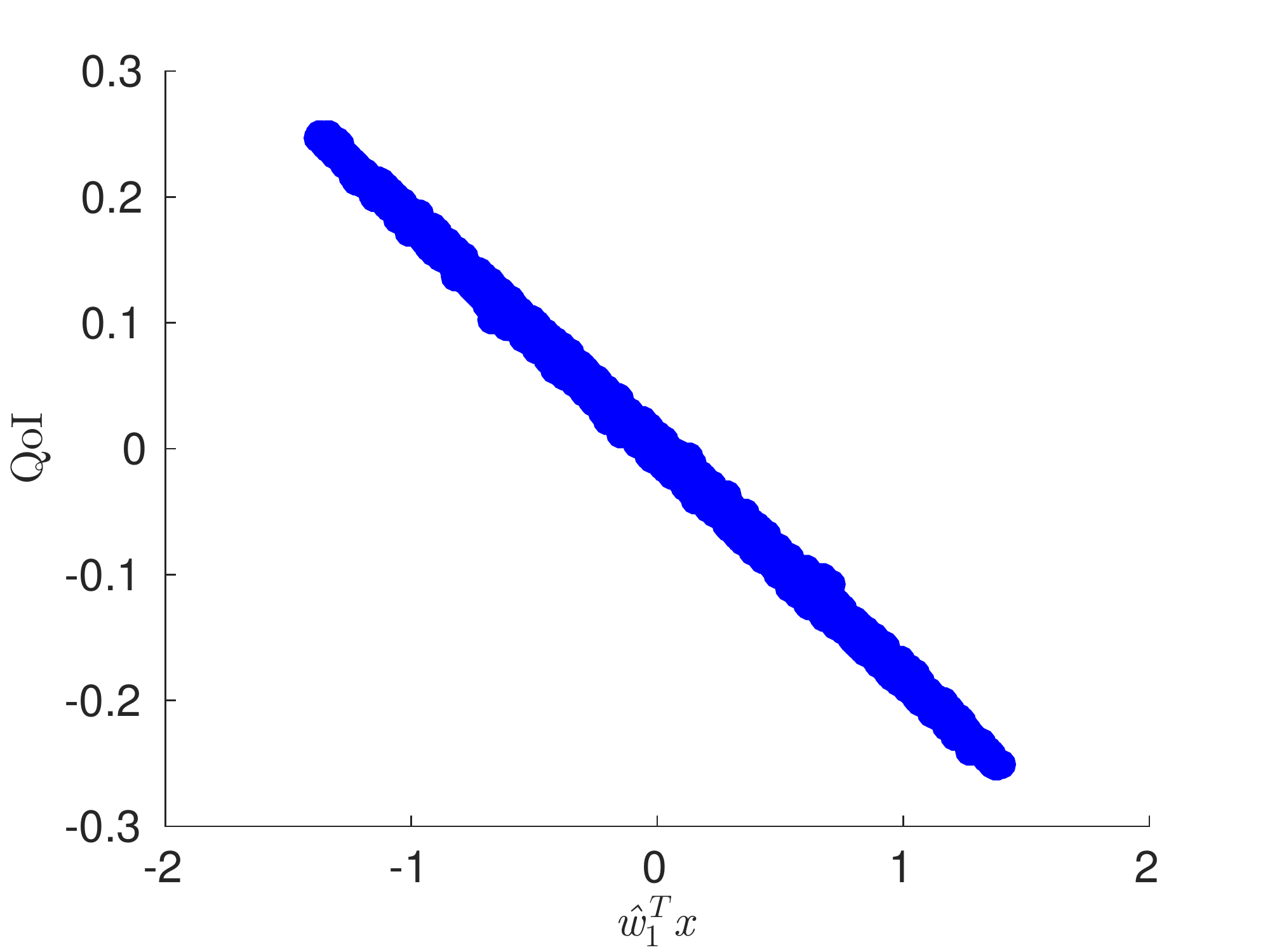}
}
\hfil
\subfloat[Drag]{
\includegraphics[width=0.45\textwidth]{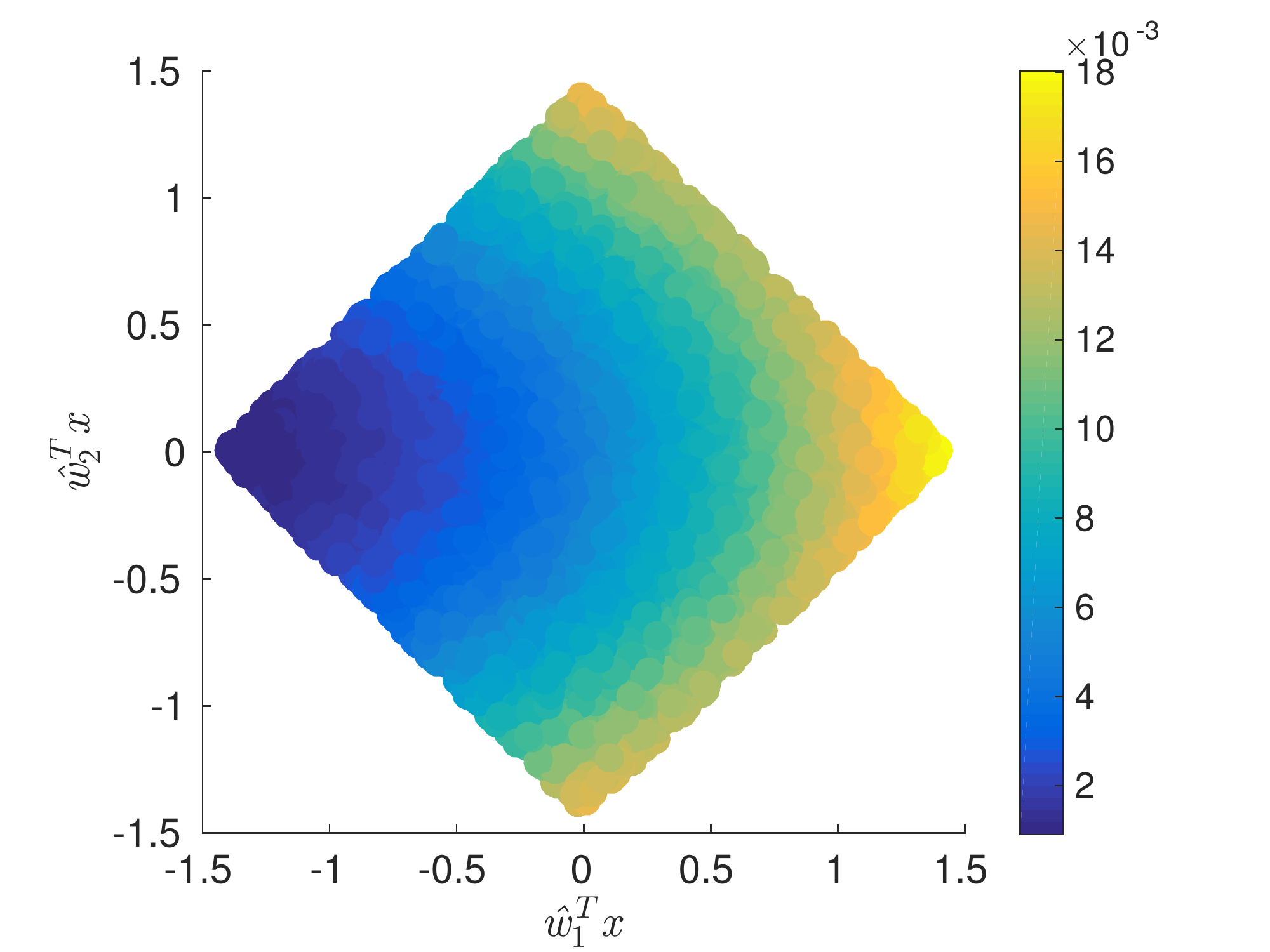}
}
\caption{One-dimensional shadow plot of lift (left) and the two-dimensional shadow plot for drag (right) using the CST shape parameterization with the first eigenvector or first and second eigenvectors from \eqref{eq:Heigs}, respectively.}
\label{CSTsum}
\end{figure}

Figure \ref{CSTsum} shows (i) the one-dimensional shadow plot for lift and (ii) the two-dimensional shadow plot for drag. Lift is very well approximated by a linear function of one linear combination of the 10 CST parameters, where the weights of the linear combination are the components of the eigenvector in Figure \ref{fig:cst-lift-evecs}. This suggests that lift can be approximated as in \eqref{eq:gfunc} using the first eigenvector from \eqref{eq:Heigs} and a linear function for $g$. Drag is very well approximated by a function of two linear combinations of the 10 CST inputs; the linear combination weights are the components of the eigenvectors in Figure \ref{fig:cst-drag-evecs}. In this case, drag can be approximated as in \eqref{eq:gfunc} using the first two eigenvectors from \eqref{eq:Heigs} and a quadratic function for $g$. The fact that the relationship between the active variables and drag appears quadratic is independent from the quadratic model used to estimate the active subspaces. We emphasize that the drag values in Figure \ref{CSTsum} are computed with the SU2 flow solver---not the quadratic response surface used to estimate the eigenvectors from \eqref{eq:Heigs}. 

Since the eigenvectors' nonzero components correspond to the leading coefficients in the CST series \eqref{eq:cstseriest}, the active subspaces are directly related to the circular leading edge radii; see Figure \ref{fig:round}. According to the eigenvector defining lift's one-dimensional active subspace (Figure \ref{CSTvec}), the subspace relates to shape perturbations from equal and opposite changes in the upper and lower surface radii. That is, positive perturbations to the active variable increase the upper surface radius and decrease the lower surface radius by the same magnitude. This perturbation agrees with the intuition that changes in camber modify lift the most on average. Additionally, the first eigenvector for drag indicates exactly equal perturbations to the leading edge radii (i.e., positive perturbations to the first active variable increase the upper and lower surface radii equally). The resulting perturbation agrees with the intuition that changes in thickness modify drag the most on average. However, the PARSEC parameterization subspaces both indicated non-zero weights corresponding to trailing edge perturbations of the airfoil shape which were omitted from the CST parameterization. Therefore, repeating the CST study to incorporate trailing edge perturbations may reveal subspaces different from those defined by the eigenvectors in Figure~\ref{CSTvec}.

\subsection{Convergence of active subspaces with number of samples}
\begin{figure}[H]
\centering
\subfloat[Lift]{
\includegraphics[width=0.45\textwidth]{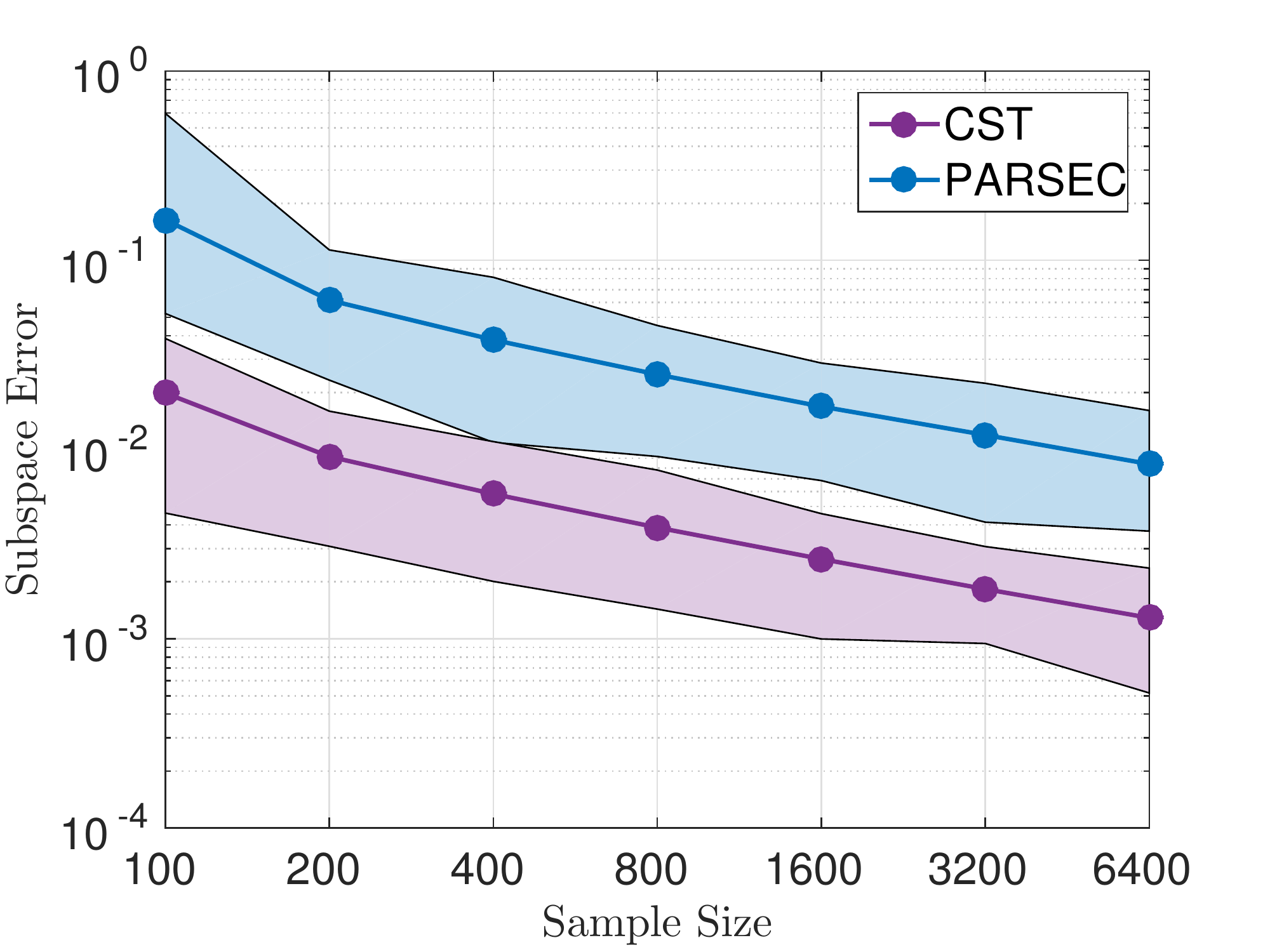}
}
\hfil
\subfloat[Drag]{
\includegraphics[width=0.45\textwidth]{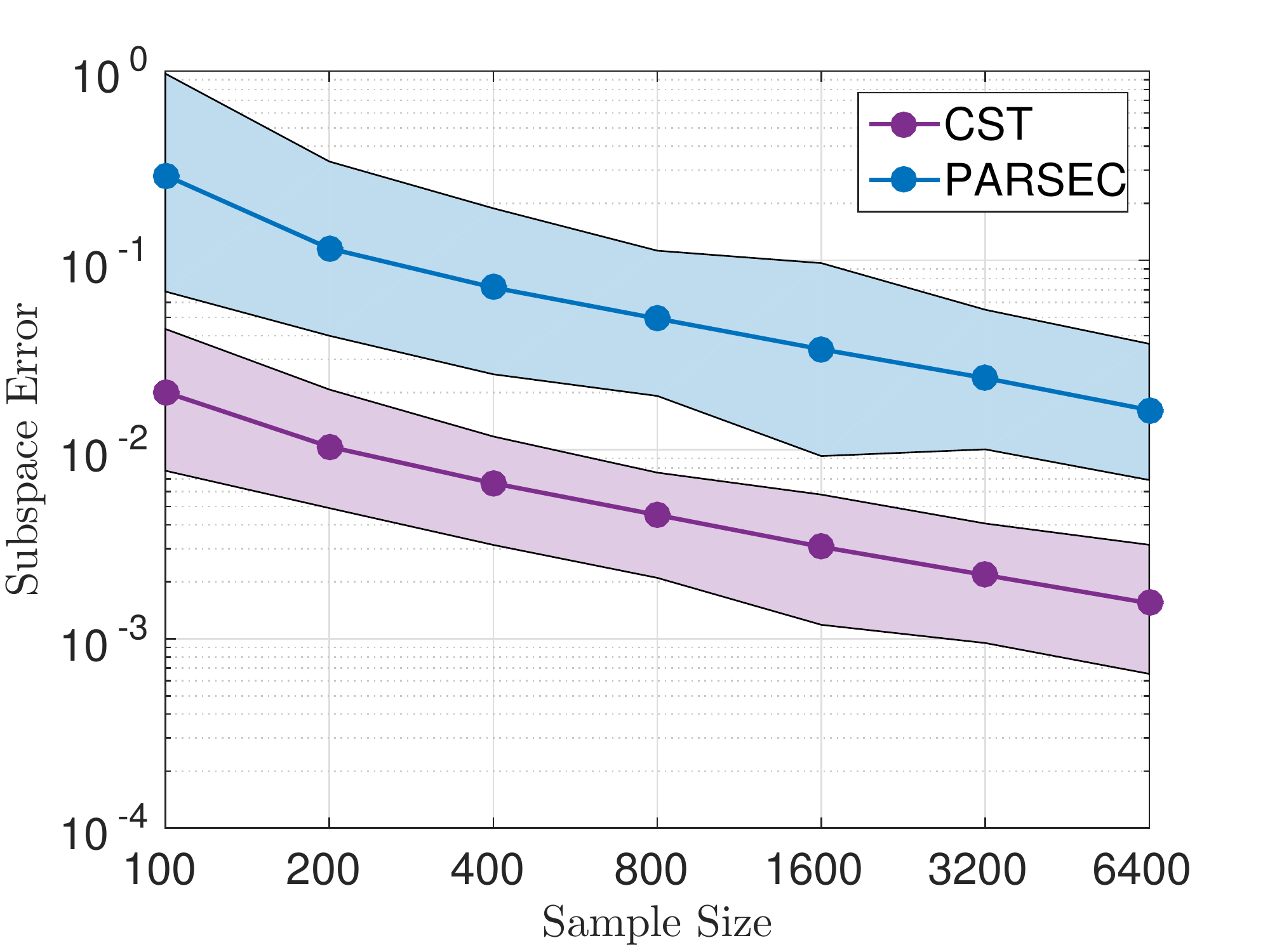}
}
\caption{Subspace error estimates (means and ranges) from \eqref{eq:booterror} for active subspaces estimated with the quadratic model-based approach (section \ref{ssec:quad}) for one-dimensional lift subspaces (left) and two-dimensional drag subspaces (right) using indicated shape parameterization for various sample sizes.}
\label{conv}
\end{figure}

The conclusions about the input/output relationships drawn from the active subspace analysis are dependent upon the number $N$ of samples used to estimate the active subspaces (in this case, via the least-squares fit global quadratic model). We study the convergence of the subspace errors as $N$ increases. Figure \ref{conv} shows convergence for both parameterizations using doubled sample sizes from 100 up to 6400.
  
The intervals in Figure \ref{conv} are the maximum and minimum from the bootstrap routine in \ref{list:bootrout}. The point values were computed as averages over the bootstrap samples. For increasing sample size, $N$, there is approximately $1/\sqrt{N}$ convergence using the global quadratic approximation of the subspaces. The subspace errors for all doubled sample sizes were consistently less with CST compared to PARSEC. Therefore, for the same number of samples, the CST parameterization resulted in more accurate subspace estimates, on average. 

\subsection{Multi-objective design interpretation}
\begin{figure}[H]
	\centering
	\subfloat[Drag]{
		\includegraphics[width=0.45\textwidth]{./CD/CST/summary2}
	}
	\hfil
	\subfloat[Lift and Drag]{
		\includegraphics[width=0.45\textwidth]{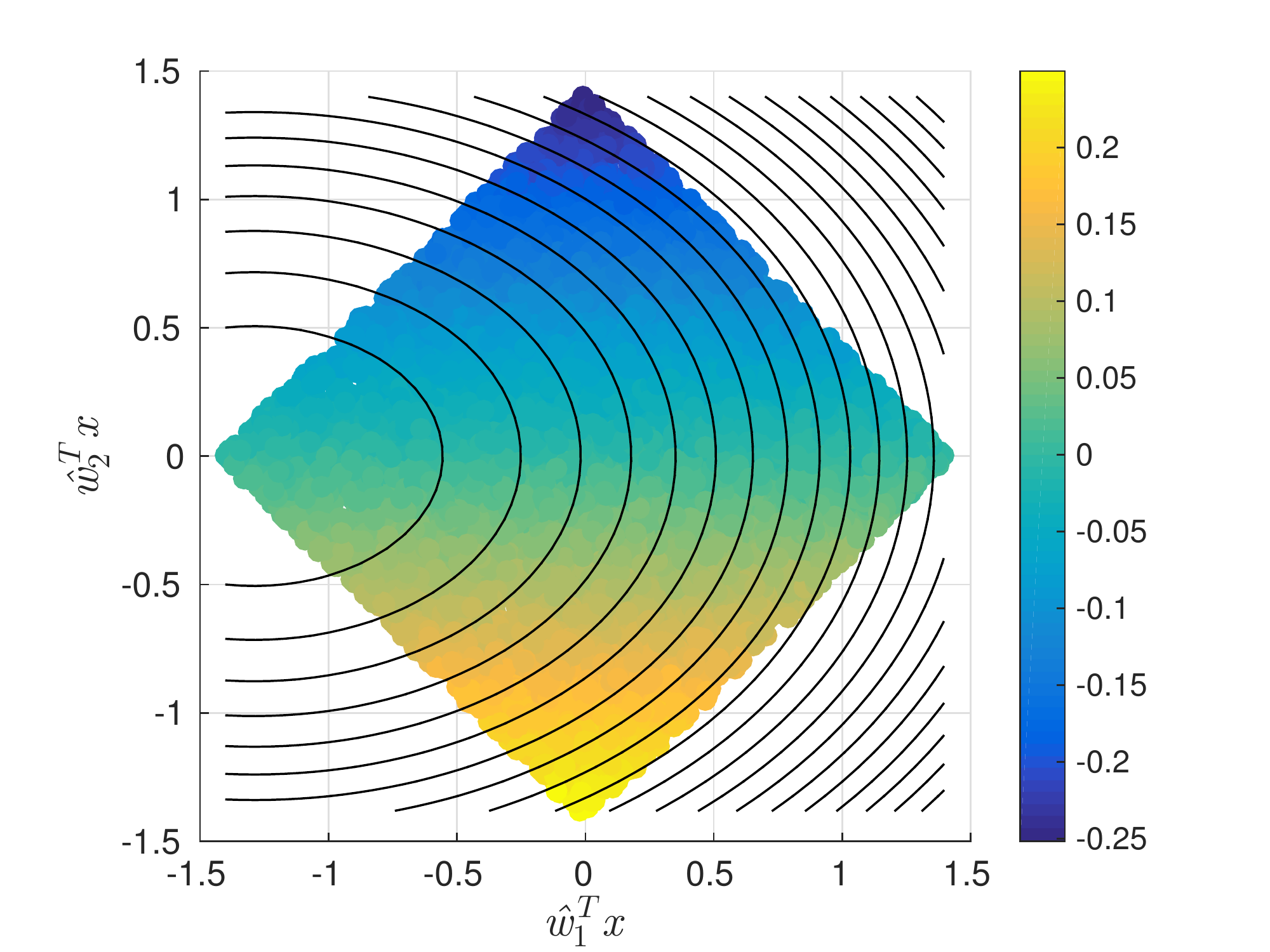}
	}
	\caption{Two-dimensional shadow plots of drag (left) and lift (right) using the first two eigenvectors of \eqref{eq:Heigs} associated with drag. The right figure also shows contours of drag overlaying the shadow plot for lift.}
	\label{multobj}
\end{figure}

For the purpose of design, we want to exploit the clear univariate and bivariate relationships between the active variables (linear combinations of the parameters) and the outputs present using the CST parameterization (see Figure \ref{CSTsum}). Leveraging the dimension reduction, we can visualize solutions to a multi-objective optimization problem. Assume the goal is to minimize drag and maximize lift. Recall that the second eigenvector for drag was the same as the eigenvector for lift; see Figure \ref{CSTvec}. Therefore, a two-dimensional shadow plot of lift versus the linear combinations defined by the drag eigenvectors shows variation in lift along the second active variable only; see Figure \ref{multobj}.

Figure \ref{multobj} shows the two-dimensional shadow plot of lift using the two eigenvectors that define drag's two-dimensional active subspace. Since lift's one-dimensional active subspace is also defined by drag's second eigenvector, the shadow plot shows that lift varies linearly along the second drag active variable. The contours overlaying the lift shadow plot show drag as a function of its active variables. The perceived relationship between the lift shadow and the drag contours immediately reveals the Pareto front for the multi-objective design problem\cite{Vanderplaats2007}. In particular, the lower left edge of the shadow plot contains the set of drag active variables $\vy=(\hvw_1^T\vx,\hvw_2^T\vx)$, constrained by the parameter space from Table \ref{CSTtable}, that either (i) maximize lift for fixed drag or (ii) minimize drag for fixed lift. Let $P_{\vy}$ denote this set of active variables,
\begin{equation}
\label{eq:paretoy}
P_{\vy} \;=\; 
\left\{\,
\vy\in\mathbb{R}^2 \,\mid\, \vy=\gamma\,\bmat{y_{1,\text{min}}\\ 0} + (1-\gamma)\,\bmat{0 \\ y_{2,\text{min}}}, \; \gamma\in[0,1]
\,\right\},
\end{equation}
where, for drag eigenvectors $\hvw_1$ and $\hvw_2$,
\begin{equation}
y_{1,\text{min}}=\mini{\vx\in[-1,1]^{10}}\,\hvw_1^T\vx, \quad
y_{2,\text{min}}=\mini{\vx\in[-1,1]^{10}}\,\hvw_2^T\vx.
\end{equation}
The hypercube parameter space $[-1,1]^{10}$ is a shifted and scaled version of the hyperrectangle defined by the ranges in Table \ref{CSTtable}. 

\begin{figure}[H]
\centering
\includegraphics[width=0.65\textwidth]{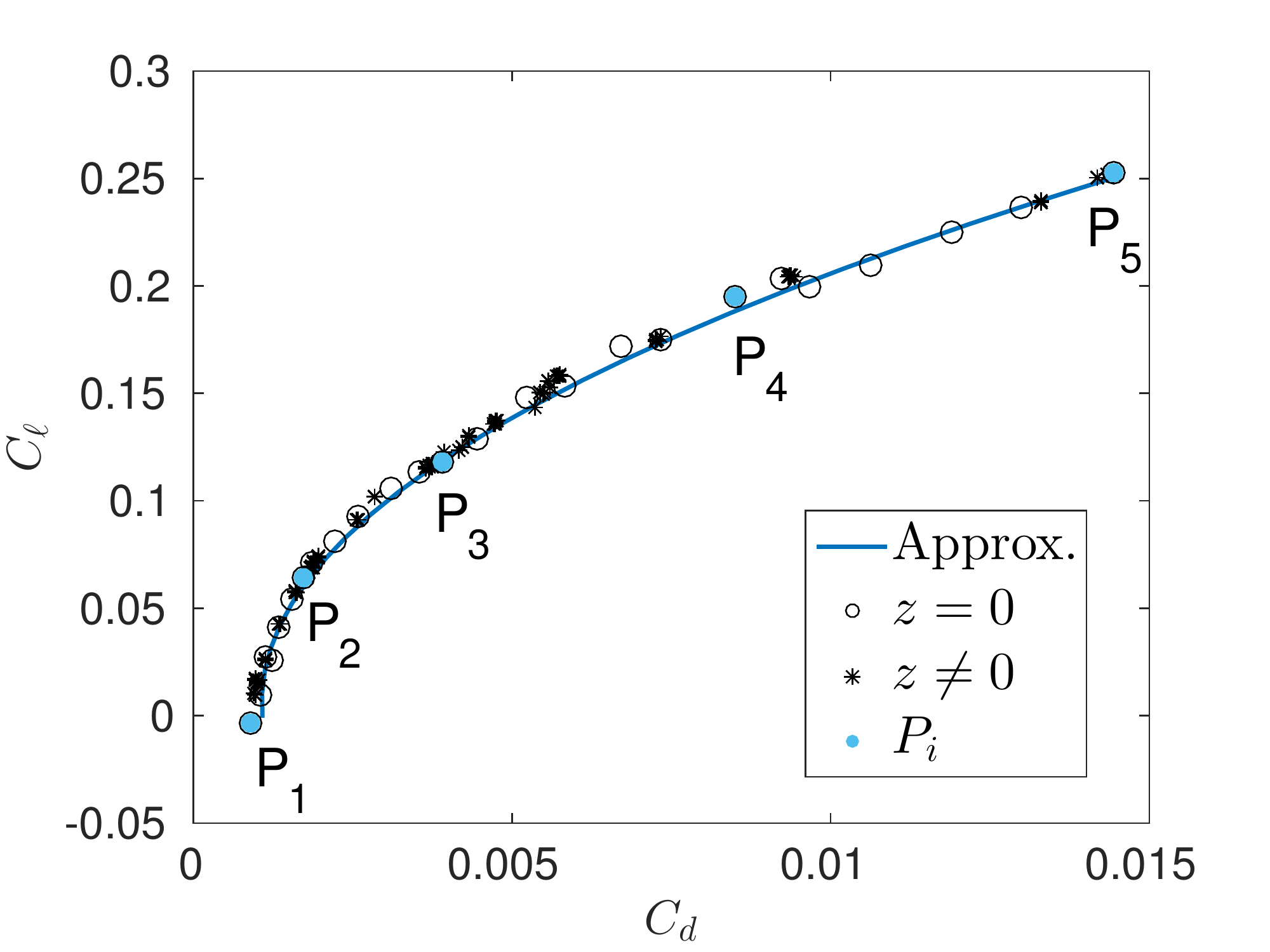}
\caption{The blue line is the lift/drag Pareto front estimated with low-dimensional active subspace-exploiting response surfaces. The circles and asterisks show lift/drag combinations computed with SU2 at points whose active variables are in the set $P_{\vy}$ from \eqref{eq:paretoy}. The circles correspond to a choice of $\vz=0$ and the asterisks are random $\vz$ such that the resulting $\vx$ is in the bounded parameter space. The $P_i$ for $i=1,\dots,5$ correspond to $z=0$ flow visualizations shown in Figure \ref{pareto_vis}.}
\label{pareto}
\end{figure}

We can visualize the lift/drag Pareto front using response surfaces for lift and drag as functions of the CST parameters. The shadow plots in Figure \ref{CSTsum} suggest simple response surfaces of the form \eqref{eq:gfunc} for each: (i) lift is well approximated by a linear function of its first active variable and (ii) drag is well approximated by a quadratic function of its first two active variables. The blue line in Figure \ref{pareto} shows lift versus drag---computed with response surfaces---for all active variables in $P_{\vy}$ from \eqref{eq:paretoy}. 

Note that for any value of the active variables $\vy$, there are infinitely many $\vx$ in the parameter space---and consequently infinitely many CST shapes---that map to the same $\vy$. Each $\vx$ that maps to a fixed $\vy$ has different values of the inactive variables $\vz$; see \eqref{paramsub}. However, by construction, we expect that perturbing the inactive variables does not change the outputs (lift and drag) very much. We verify this insensitivity by choosing several values for $\vz$ for each $\vy$ such that the resulting $\vx$ are in the parameter space defined by the ranges in Table \ref{CSTtable}. Each $\vy,\vz$ combination produces a unique airfoil shape. In Figure \ref{pareto}, several SU2-computed lift/drag combinations are shown; the CST shape for each combination has a parameter vector $\vx$ whose active subspace coordinates $\vy=\hmW_1^T\vx\in\mathbb{R}^2$ are in $P_{\vy}$ from \eqref{eq:paretoy}. The circles correspond to a particular choice of $\vz=0$, and the asterisks correspond to random $\vz$ values such that $\vx=\hmW_1\vy + \hmW_2\vz$ is in the (scaled) parameter space. Note that all SU2-computed lift/drag combinations are remarkably close to the response surface-generated Pareto front. This verifies that (i) the inactive variables for lift and drag are indeed relatively inactive and (ii) the Pareto front estimated with quadratic model-based active subspaces and active subspace-exploiting low-dimensional response surfaces for lift and drag is accurate. We emphasize that this conclusion depends on several assumptions including the density $\rho$ derived from Table \ref{CSTtable} and the angle-of-attack and NIST conditions for the flow field simulations. If any of these assumptions change, then the numerical study should be repeated.

\begin{figure}[H]
\centering
\includegraphics[width=0.65\textwidth]{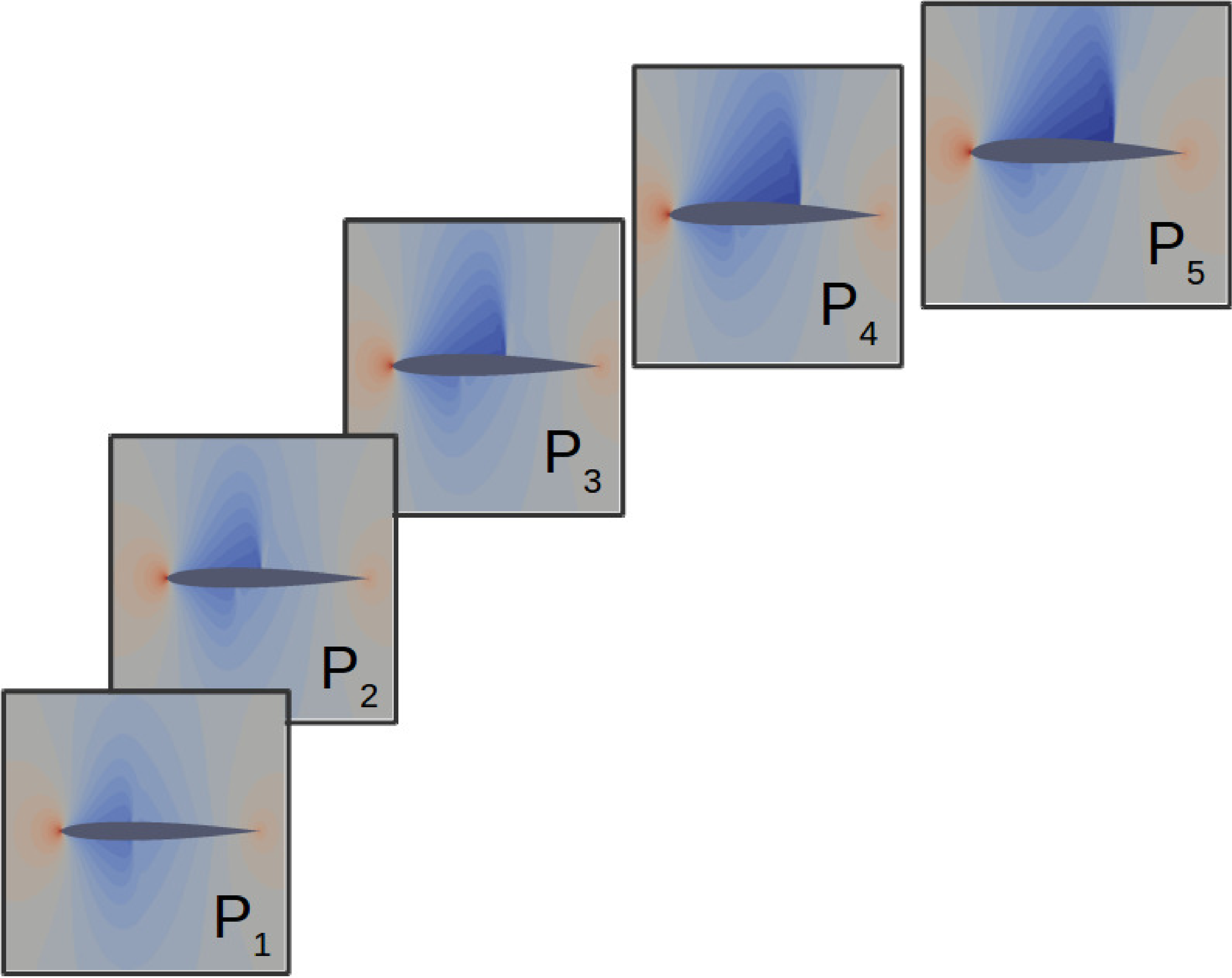}
\caption{Five Paraview flow visualizations computed with SU2 at points whose active variables are in the set $P_{\vy}$ from \eqref{eq:paretoy} also shown in Figure \ref{pareto}. These particular points correspond to a choice of $\vz=0$. }
\label{pareto_vis}
\end{figure}

In general, two nearby points on the Pareto front do not have similar input parameters or airfoil shapes\cite{Vanderplaats2007}. There is no notion of continuity among the shapes whose lift/drag combinations define the Pareto front. Standard surrogate-based approaches\cite{Simpson2001} offer no means to reveal continuous perturbations to airfoil shapes that remain on the Pareto front. With active subspaces, we can define a conditional Pareto front as a function of the active variables by averaging over the inactive variables $\vz$ or choosing a fixed $\vz$ of interest. Then airfoil shapes are continuous with changes in the active variables. This continuity in the shape is emphasized in Figure \ref{pareto_vis} as five flow visualizations corresponding to the $P_i$ points in Figure \ref{pareto}.

For example, taking $\vz=0$, the values $\vy$ in the Pareto subset $P_{\vy}$ from \eqref{eq:paretoy}---i.e., active variables along the lower left edge of the shadow plot in Figure \ref{multobj}---produce continuous changes to the airfoil shape parameters and consequently continuous changes in the airfoil shape. We expect this continuity, when enabled by the active subspaces, would be very helpful for airfoil designers. 


\section{Summary and conclusions}

We have examined two transonic airfoil shape parameterizations. PARSEC uses a set of parameterized equality constraints to define the coefficients of a series that defines the upper and lower airfoil surfaces; each parameter affects a physical characteristic of the surfaces. CST's parameters define the series' coefficients directly, constrained to satisfy leading and trailing edge conditions. For each parameterization, we estimate active subspaces for airfoil performance metrics lift and drag---computed with CFD---as functions of the airfoil parameters. To estimate the active subspaces, we employ a global least-squares-fit quadratic approximation. The active subspaces derived from this analysis reveal one-dimensional structure in lift and two-dimensional structure in drag for each shape parameterization---the 11-parameter PARSEC and the 10-parameter CST. The vectors defining the active subspaces give insight into the important shape  parameters for each performance metric, and these insights agree with first principles aerodynamics in both parameterizations. Moreover, the low-dimensional active subspaces enable simple visualization tools called shadow plots that aid in response surface construction. Additionally, the visualizations and low-dimensional active subspaces reveal insights into the Pareto front of a multi-objective design problem using both lift and drag as objective functions. For example, we can describe continuous perturbations to the airfoil that move smoothly along the Pareto front.

\pagebreak

\bibliographystyle{plain}
\bibliography{jab_database}
\pagebreak
\appendix
\section{Appendix}
\label{sec:app1}
The appendix details the PARSEC constraints for a single shape function and an alternative leading edge radius relationship proof.
\\
\\
Constraints:
\begin{enumerate}

\item For a single interior interpolation point $ \ell_{int} \in (0,1)$:
$$
s(\ell_{int}) = \boldsymbol{\phi}(\ell_{int})^T\boldsymbol{a}
$$
\item Trailing edge point interpolation $\ell=1$:
$$
s(1) = \boldsymbol{\phi}(1)^T\boldsymbol{a} \implies \boldsymbol{1}^T\boldsymbol{a}
$$
\item Interior interpolation point slope of zero $\ell_{int} \in (0,1)$:
$$
0=\frac{d s(\ell_{int})}{d\ell} = \left(\frac{d }{d\ell}\left[\boldsymbol{\phi}(\ell_{int})\right]\right)^T \boldsymbol{a}
$$
\item Trailing edge slope $\ell=1$:
$$
\frac{d s(1)}{d\ell} = \left(\frac{d }{d\ell}\left[\boldsymbol{\phi}(1)\right]\right)^T \boldsymbol{a}
$$
\item Interior point inflection at $\ell_{int} \in (0,1)$:
$$
\frac{d^2 s(\ell_{int})}{d\ell^2} = \left(\frac{d^2}{d\ell^2}\left[\boldsymbol{\phi}(\ell_{int})\right]\right)^T\boldsymbol{a}
$$

\item Intermediate scale circular approximation $0<\epsilon \ll 1$, $\ell \in (0,\epsilon)$:
$$
\sqrt{2\epsilon} = \boldsymbol{e}_1^T\boldsymbol{a}  \quad \text{where} \quad \boldsymbol{e}_1 = (1,\,\, 0,\,\, 0,\,\, 0,\,\, 0,\,\, 0)^T
$$
\end{enumerate}
For the six constraints, we construct a linear system to represent the separate shape function equality constraints 
$$
\mM \boldsymbol{a} = \boldsymbol{p} \quad \text{where} \quad \boldsymbol{p} = \left(s(\ell_{int}), s_i(1), 0, \frac{d s_{int}(1)}{d\ell}, \frac{d^2 s_{int}(\ell_{int})}{d\ell^2},\sqrt{2\epsilon}\right)^T,
$$
and invertible
$$
\mM = \left[ \begin{matrix}
\sqrt{\ell_{int}} & \left(\sqrt{\ell_{int}}\right)^3 & \left(\sqrt{\ell_{int}}\right)^5 & \left(\sqrt{\ell_{int}}\right)^7 & \left(\sqrt{\ell_{int}}\right)^9 & \left(\sqrt{\ell_{int}}\right)^{11}\\[5pt] 
1 & 1 & 1 & 1 & 1 & 1 \\[5pt]
\frac{1}{2}\left(\sqrt{\ell_{int}}\right)^{-1} & \frac{3}{2}\sqrt{\ell_{int}} & \frac{5}{2}\left(\sqrt{\ell_{int}}\right)^{3}& \frac{7}{2}\left(\sqrt{\ell_{int}}\right)^{5} & \frac{9}{2}\left(\sqrt{\ell_{int}}\right)^{7} & \frac{11}{2}\left(\sqrt{\ell_{int}}\right)^{9}\\[5pt]
\frac{1}{2} & \frac{3}{2} & \frac{5}{2} & \frac{7}{2} & \frac{9}{2} & \frac{11}{2}\\[5pt]
-\frac{1}{4}\left(\sqrt{\ell_{int}}\right)^{-3} & \frac{3}{4}\sqrt{\ell_{int}}^{-1} & \frac{15}{4}\left(\sqrt{\ell_{int}}\right)& \frac{35}{4}\left(\sqrt{\ell_{int}}\right)^{3} & \frac{63}{4}\left(\sqrt{\ell_{int}}\right)^{5} & \frac{99}{4}\left(\sqrt{\ell_{int}}\right)^{7}\\[5pt]
1 & 0 & 0 & 0 & 0 & 0\\[5pt]
\end{matrix}\right].
$$


Proof of leading coefficient and leading edge radius dependency:
\begin{proof}
Determining the appropriate coefficient for the sixth constraint requires us to consider intermediate scale $t$. Consider a  change of variables, $t = \sqrt{\ell}$, and recall the odd-power polynomial series expansion as the basis for PARSEC with six terms and CST with arbitrary basis,
$$
s(t; \boldsymbol{a}) = \sum_{j=1}^{k} a_jt^{2j-1} \quad \text{for all} \quad t \in [0,1].
$$
Additionally, construct a leading edge circular approximation with radius $\epsilon$ at a distance $\epsilon$ from the origin such that
$$
s_c^2 + (\ell-\epsilon)^2 = \epsilon^2.
$$
Substituting the intermediate scale change-of-variables and solving for $s_c^2(t)$ results in
$$
s_c^2(t) = 2\epsilon\ell-\ell^2 \implies 2\epsilon t^2 -t^4.
$$
Truncate the odd polynomial basis expansion for $k>2$ and square the expression such that
$$
s(t) = a_1t + a_2t^3 \implies s^2(t) = a_1^2t^2 + 2a_1a_2t^4+a_2^2t^6
$$
For $0 < \epsilon \ll 1$, assume $t^2 \in (0,\epsilon]$ (i.e., $\ell \in (0,\epsilon]$). Truncating powers of $t$ greater than $2$ results in the small scale approximations of the circular function and expansion,
$$
s_c(t) \approx 2\epsilon t^2 \quad \text{and} \quad s(t) \approx a_1^2 t^2.
$$
Equating the two approximations and solving for the coefficient when $t \neq 0$ implies
$$
a_1 = \pm \sqrt{2\epsilon}.
$$ 
The positive sign corresponds to the upper surface coefficient while the negative sign corresponds to the lower surface coefficient.
\end{proof}

\end{document}